\documentclass{aic}

\aicAUTHORdetails{%
  title = {Improved bounds for the Erd\H os-Rogers function}, 
  author = {W. T. Gowers and O. Janzer},
  plaintextauthor = {W. T. Gowers, O. Janzer},
  plaintexttitle = {Improved bounds for the Erdos-Rogers function}, 
  keywords = {Ramsey theory},
}   

\aicEDITORdetails{%
   year={2020},
   number={3},
   received={25 July 2018},   
   revised={11 June 2019},    
   published={28 February 2020},  
   doi={10.19086/aic.12048},      
}   

\usepackage{pdfpages}
\usepackage{graphicx}
\usepackage{mathtools}
\usepackage{verbatim}
\usepackage{amssymb}
\usepackage{bm}
\usepackage{enumitem}
\usepackage{cases}
\usepackage{latexsym}
\usepackage{amsthm}
\usepackage{amsfonts}
\usepackage{bbm}
\usepackage{etoolbox}
\usepackage[top=25mm, bottom=25mm, left=28mm, right=28mm]{geometry}
\usepackage{cite}
\usepackage{changepage}

\newtheorem{theorem}{Theorem}
\AfterEndEnvironment{theorem}{\noindent\ignorespaces}
\newtheorem{corollary}[theorem]{Corollary}
\AfterEndEnvironment{corollary}{\noindent\ignorespaces}
\newtheorem{lemma}[theorem]{Lemma}
\AfterEndEnvironment{lemma}{\noindent\ignorespaces}

\AfterEndEnvironment{proposition}{\noindent\ignorespaces}
\newtheorem*{question}{Question}
\AfterEndEnvironment{question}{\noindent\ignorespaces}

\theoremstyle{definition}
\newtheorem{definition}[theorem]{Definition}
\AfterEndEnvironment{definition}{\noindent\ignorespaces}

\theoremstyle{remark}
\newtheorem*{example}{Example}
\AfterEndEnvironment{example}{\noindent\ignorespaces}
\newtheorem*{remark}{Remark}
\AfterEndEnvironment{remark}{\noindent\ignorespaces}

\AfterEndEnvironment{notation}{\noindent\ignorespaces}

\newenvironment{proof2}{\proof[\textnormal{\textbf{Proof.}}]}{\qed}

\usepackage{chngcntr}
\usepackage{apptools}
\AtAppendix{\counterwithin{theorem}{section}}

\def\Bin{\mathop{\mathrm{Bin}}}

\begin{document}

\begin{frontmatter}[classification=text]

\title{Improved bounds for the Erd\H os-Rogers function} 

\author[tim]{W. T. Gowers\thanks{Royal Society 2010 Anniversary Research Professor, University of Cambridge}}
\author[oliver]{O. Janzer\thanks{Department of Pure Mathematics and Mathematical Statistics, University of Cambridge}}

\begin{abstract}
The Erd\H os-Rogers function $f_{s,t}$ measures how large a $K_s$-free induced subgraph there must be in a $K_t$-free graph on $n$ vertices. While good estimates for $f_{s,t}$ are known for some pairs $(s,t)$, notably when $t=s+1$, in general there are significant gaps between the best known upper and lower bounds. We improve the upper bounds when $s+2\leq t\leq 2s-1$. For each such pair we obtain for the first time a proof that $f_{s,t}\leq n^{\alpha_{s,t}+o(1)}$ with an exponent $\alpha_{s,t}<1/2$, answering a question of Dudek, Retter and R\"odl.
\end{abstract}
\end{frontmatter}

\section{Introduction}

Let $G$ be a graph with $n$ vertices that contains no $K_4$. How large a triangle-free induced subgraph must $G$ have? The standard proof of Ramsey's theorem implies that $G$ contains an independent set of size $n^{1/3}$, but can we do better? 

A simple argument shows that the answer is yes. Indeed, each vertex in $G$ has a triangle-free neighbourhood, and either there is a vertex of degree $n^{1/2}$ or one can find an independent set of size roughly $n^{1/2}$ by repeatedly choosing vertices and discarding their neighbours. 

This stronger argument still feels a little wasteful, because in the second case one finds an independent set rather than a triangle-free subgraph. Moreover, there is no obvious example that yields a matching upper bound, so it is not immediately clear whether 1/2 is the correct exponent.

The problem above is an example of a general problem that was first considered by Erd\H os and Rogers. Given positive integers $1<s<t$ and $n>2$, define $f_{s,t}(n)$ to be the minimum over all $K_t$-free graphs $G$ with $n$ vertices of the order of the largest induced $K_s$-free subgraph of $G$. We have just been discussing the function $f_{3,4}$. The function $f_{s,t}$ is known as the \emph{Erd\H os-Rogers function}. It has been studied by several authors: for a detailed survey covering many of the known results on the subject, see \cite{dudekrodlsurvey}. For a more recent exposition, see also section 3.5.2 of \cite{sudakovsurvey}.

\medskip
The first bounds were obtained by Erd\H{o}s and Rogers \cite{erdosrogers} who showed that for every $s$ there exists a positive constant $\epsilon(s)$ such that $f_{s,s+1}(n)\leq n^{1-\epsilon(s)}$. About 30 years later, Bollob\'{a}s and Hind \cite{bollobashind} improved the estimate for $\epsilon(s)$ and established the lower bound $f_{s,t}(n)\geq n^{1/(t-s+1)}$. In particular, $f_{s,s+1}(n)\geq n^{1/2}$ (by the obvious generalization of the argument for $f_{3,4}$ above).

Subsequently, Krivelevich \cite{krivelevichlower,krivelevichupper} improved these lower bounds by a small power of $\log n$ and also gave a new general upper bound, which is

\begin{equation} \label{kriv}
	f_{s,t}(n)\leq O(n^{\frac{s}{t+1}}(\log n)^{\frac{1}{s-1}}).
\end{equation}

Later, the lower bound was significantly improved by Sudakov \cite{sudakovold,sudakovnew}. He showed that if $t>s+1$, then $f_{s,t}(n)\geq \Omega(n^{a_{s,t}})$ where $a_{s,t}$ is defined recursively. In particular, when $s$ is fixed and $t\rightarrow \infty$, he obtained the bound

\begin{equation*} \label{sud}
	f_{s,t}(n)\geq \Omega(n^{\frac{s}{2t}+O(1/t^2)}).
\end{equation*}

We remark that if $t\geq 2s$ then (\ref{kriv}) is the best known upper bound, while Sudakov's lower bound is the best known for every $t>s+1$. In particular, the upper bound is roughly the square of the lower bound in the range $t\geq 2s$.

\vspace{2mm}

Recently, there has been quite a lot of progress on the case $t=s+1$. First, Dudek and R\"{o}dl \cite{dudekrodl} showed that $f_{s,s+1}(n)\leq O(n^{2/3})$. Then Wolfovitz \cite{wolfovitz} proved that for sufficiently large $n$ we have $f_{3,4}(n)\leq n^{1/2}(\log n)^{120}$, yielding the slightly surprising fact that the exponent 1/2 is indeed the right one in that case. Finally, Dudek, Retter and R\"{o}dl \cite{dudekretterrodl}, generalizing Wolfovitz's construction, showed that for any $s\geq 3$ there exist constants $c_1$ and $c_2$ such that

\begin{equation*}
	f_{s,s+1}(n)\leq c_1n^{1/2}(\log n)^{c_2}
\end{equation*}

\noindent so the exponent 1/2 is correct for all $f_{s,s+1}$. However, the problem of finding the correct exponent of $n$ for general $s,t$ remains open.

\vspace{2mm}

A particularly important case is when $t=s+2$ since $f_{s,t}(n)\leq f_{s,s+2}(n)$ for any $t\geq s+2$. Sudakov's lower bound gives $f_{s,s+2}(n)= \Omega(n^{\beta_s})$ where $\beta_s=1/2-\frac{1}{6s-6}$. Dudek, Retter and R\"{o}dl in \cite{dudekretterrodl} showed that for any $s\geq 4$ there exists a constant $c$ depending only on $s$ such that

\begin{equation*}
	f_{s,s+2}(n)\leq cn^{1/2}.
\end{equation*}

\noindent Note that the exponent 1/2 follows from the bound for $f_{s,s+1}$, so this improves it by removing the log factor. Having established this, Dudek, Retter and R\"odl asked the following question.

\begin{question}
	Does there exist $s\geq 3$ such that $f_{s,s+2}(n)=o(n^{1/2})$?
\end{question}

Another central open problem in the area is the following question of Erd\H{o}s \cite{erdosquestion}.

\begin{question}
	Is it true that 
	\begin{equation} \label{erdoseqn}
		\lim_{n\rightarrow \infty} \frac{f_{s+1,t}(n)}{f_{s,t}(n)}=\infty 
	\end{equation}
for every  $t>s+1$?
\end{question}

\noindent The answer has been shown to be yes when $t=s+2\geq 6$ and when $(s,t)$ is one of the pairs $(2,4)$, $(2,5)$, $(2,6)$, $(2,7)$, $(2,8)$ or $(3,6)$.

\subsection{Our results}

In this paper, we prove that the answer to the first question above is yes. We also establish (\ref{erdoseqn}) for the families of pairs $t=s+3\geq 7$ and $t=s+2\geq 5$. We obtain these results by proving a significant improvement for the upper bound on $f_{s,t}$ when $s+2\leq t\leq 2s-1$. The previous best upper bound for these parameters appeared in \cite{dudekretterrodl} and was $f_{s,t}(n)\leq cn^{1/2}$ (except for the pair $s=3,t=5$, where this bound was not established). We do not just obtain bounds of the form $o(n^{1/2})$, but we improve the exponents throughout the range. Our construction is probabilistic, and has some similarities to the constructions that established the previous best upper bounds. However, an important difference is that we do not make use of algebraic objects such as projective planes.

\vspace{3mm}

To state the bound that comes out of our argument takes a small amount of preparation. Let $s\geq 3$ and $s+2\leq t\leq 2s-1$. Call $(s,t)$ \emph{regular} if $s\geq 11$ and $s+3\leq t\leq 2s-4$ or if $(s,t)\in \{(10,14),(10,15)\}$ and call it \emph{exceptional} otherwise. Let 
\begin{equation*}
\alpha=\alpha_{s,t}=\begin{cases}
\alpha(1)=\frac{(s-2)(t-s)(t+s-1)+2t-2s}{(2s-3)(t-s)(t+s-1)-2s+4}, & \text{ if } (s,t) \text{ is regular} \\
\alpha(2)=\frac{(s-2)(t-s)(s-1)+s-1}{(2s-3)(t-s)(s-1)+2s-t}, & \text{ if } (s,t) \text{ is exceptional}
\end{cases}
\end{equation*}

\noindent We will prove the following theorem.

\begin{theorem} \label{ERbound}
	For any $s\geq 3,s+2\leq t\leq 2s-1$, there exists some constant $c=c(s,t)$ such that \begin{equation*}
	f_{s,t}(n)\leq n^{\alpha}(\log n)^c.
	\end{equation*}
\end{theorem}
It is not hard to check that $\alpha<1/2$ for all pairs $(s,t)$ in the given range. Thus, as mentioned above, we obtain a strong answer to the question of Dudek, Retter and R\"{o}dl.

\begin{corollary}
	For every $s\geq 3$, we have $f_{s,s+2}(n)=o(n^{1/2})$.
\end{corollary}

\noindent The simplest case where our result is new is the case $s=3,t=5$. There we obtain an upper bound of $n^{6/13}(\log n)^c$. For comparison, Sudakov's lower bound is $cn^{5/12}$.

Since the exponent when $t=s+1$ is 1/2, our result also implies a positive answer to the question of Erd\H os in the following family of cases.

\begin{corollary}
	\begin{equation*}
		\lim_{n\rightarrow \infty} \frac{f_{s+1,s+2}(n)}{f_{s,s+2}(n)}\rightarrow \infty
	\end{equation*}
That is, (\ref{erdoseqn}) holds for $t=s+2\geq 5$.
\end{corollary}

\vspace{3mm}

\noindent If $t=s+3$, then \begin{equation*}
	\alpha=\begin{cases}
	\frac{3s^2-3s-3}{6s^2-4s-7}, & \text{ if } s\geq 11 \\
	\frac{3s^2-8s+5}{6s^2-14s+6}, & \text{ if } 4\leq s\leq 10
	\end{cases}
\end{equation*}

\noindent Comparing this with Sudakov's lower bound $f_{s+1,s+3}(n)\geq \Omega(n^{\beta_{s+1}})$, where $\beta_{s+1}=\frac{3s-1}{6s}$, we get the following additional result.

\begin{corollary}
	\begin{equation*}
	\lim_{n\rightarrow \infty} \frac{f_{s+1,s+3}(n)}{f_{s,s+3}(n)}\rightarrow \infty
	\end{equation*}
That is, (\ref{erdoseqn}) holds for $t=s+3\geq 7$.
\end{corollary}

In the following table, we compare the exponent of $n$ in the best known lower bound with that in our new upper bound (both rounded to three decimal places).


\begin{center}
\begin{tabular} {|c||c|c|}
\hline
& Our new upper bound & Best known lower bound\\
\hline
\hline
$s=3, t=5$ & 0.462 & 0.417\\
\hline
$s=4,t=6$ & 0.467 & 0.444\\
\hline
$s=4,t=7$ & 0.457 & 0.375\\
\hline
$s=5,t=7$ & 0.475 & 0.458\\
\hline
$s=5,t=8$ & 0.465 & 0.404\\
\hline
$s=5,t=9$ & 0.460 & 0.351\\
\hline
\end{tabular}
\end{center}

\noindent In the case $t=s+2$, our bound is $f_{s,s+2}(n)\leq n^{\alpha+o(1)}$ for $\alpha=1/2-\frac{s-2}{8s^2-18s+8}\approx 1/2-\frac{1}{8s}$ while Sudakov's lower bound is $f_{s,s+2}(n)\geq n^{\beta+o(1)}$ for $\beta=1/2-\frac{1}{6s-6}\approx 1/2-\frac{1}{6s}$. It would be very interesting to know whether either of these two estimates reflects the true asymptotics of $f_{s,s+2}$. It would be particularly interesting to know whether either of the exponents $5/12$ or $6/13$ is the correct one for $(s,t)=(3,5)$. We have made some effort to optimize our construction, whereas there appear to be places where Sudakov's argument is potentially throwing information away, so our guess is that $6/13$ is correct, but this guess is very tentative and could easily turn out to be wrong.

\subsection{An overview of the argument}

We will now sketch the key steps in our argument. For simplicity, we will focus on the $s=3,t=5$ case. Then, as mentioned above, Theorem \ref{ERbound} says that $f_{3,5}(n)\leq n^{6/13}(\log n)^c$. That is, we construct a $K_5$-free graph $G$ in which every subset of size roughly $n^{6/13}$ induces a triangle.

The basic idea is very simple. We are looking for a graph that contains ``triangles everywhere" but does not contain any $K_5$s. The obvious way to create a large number of triangles without creating $K_5$s is to take a complete tripartite graph. Of course, this on its own does nothing, since a complete tripartite graph has a huge independent set, but we can use it as a building block by taking a union of many complete tripartite graphs. In previous constructions, such as Wolfovitz's graph that gives an upper bound for $f_{3,4}(n)$, the vertex sets of these tripartite graphs are chosen algebraically -- in Wolfovitz's case they are the lines of a projective plane. The main difference in our approach is that we simply choose them at random, where the number we choose and the size of each one are parameters that we optimize at the end of the argument. This creates difficulties that are not present in the earlier approaches, but in the end allows us to prove stronger bounds. 

Thus, we begin by taking a graph $G_0$, which is a union of roughly $n^{9/13}$ complete tripartite graphs with parts having size roughly $n^{6/13}$ each, these parts being randomly chosen subsets of $V(G_0)$. It is not hard to prove that $G_0$ contains a triangle in every set of vertices of size roughly~$n^{6/13}$. 

However, $G_0$ also contains many $K_5$s, so we have to delete some edges. It is here that the proof becomes less simple: while random constructions followed by edge deletions are very standard, in this case we need rather delicate arguments in order to prove that it can be done without removing all the triangles from a set of size around $n^{6/13}$.

First, let us check that every set of size roughly $n^{6/13}$ does indeed induce a triangle in $G_0$. Let $A$ be a subset of $V(G_0)$ of size $n^{6/13}$. A given tripartite copy will intersect $A$ in at least 3 vertices with probability roughly $n^{-3/13}$. Thus, as we place $n^{9/13}$ tripartites, the expected number of those tripartites that give a triangle in $A$ is roughly $n^{6/13}$. Hence, by the Chernoff bound, the probability that $A$ does not contain a triangle is roughly $e^{-n^{6/13}}$. But the number of subsets of $V(G_0)$ of size $n^{6/13}$ is \textit{very} roughly $n^{n^{6/13}}$. Modifying the parameters by $\log n$ factors suitably, a union bound shows that almost surely every subset $A$ of size roughly $n^{6/13}$ will contain a triangle. In fact, a slightly more careful examination of this argument reveals that almost surely every such subset will contain at least $n^{6/13}$ triangles, each coming from a single tripartite graph such that the tripartites corresponding to different triangles are all distinct.

Now let us specify which edges get deleted. We shall delete them in two stages. The first stage consists of what we call \emph{Type 1} deletions. Given any two of our random tripartite graphs, with vertex sets $A=A_0\cup A_1\cup A_2$ and $B=B_0\cup B_1\cup B_2$, we remove all edges $xy$ such that $x,y\in A\cap B$. We do not insist that $xy$ is an edge of both tripartite graphs: if, for example, $x,y\in A_0$, $x\in B_0$ and $y\in B_1$, then the edge $xy$ will be removed. Let $G_1$ be the resulting graph when all such edges have been deleted. The reason for these deletions is that each of our tripartite graphs contains many copies of $K_{3,1,1}$, which are somewhat ``dangerous" for us, since all it takes to convert a $K_{3,1,1}$ into a $K_5$ is the addition of a further triangle. If we do not do Type 1 deletions, then we will obtain $K_5$s in this way too frequently, with the result that most edges in the graph are contained in a $K_5$. Indeed, the expected number of edges in $G_0$ is roughly $n^{9/13}(n^{6/13})^2=n^{21/13}$ and the expected number of $K_5$s of the above form is roughly $n^5(n^{9/13})^2(n^{-7/13})^8=n^{27/13}$.

Type 1 deletion is feasible in the sense that it destroys only a small proportion of the edges of $G_0$. That is because it is significantly less likely for a pair of vertices to be contained in two tripartite copies than for it to be contained in one tripartite copy.

Thanks to Type 1 deletions, it has become ``difficult" for $K_5$s to appear in $G_1$, since now none of our random tripartite graphs can intersect a $K_5$ in more than 3 vertices. Indeed, if one of them intersects a $K_5$ in say 4 vertices, then there exist two of those vertices between which this tripartite does not provide an edge, and if one of the other tripartites gives an edge in $G_0$ between those two vertices, that edge is deleted.

Thus, it is easy to check that if a $K_5$ appears in $G_1$, then it has to do so in one of the following ways.
\begin{enumerate}[label=(\roman*)]
\item All 10 edges of the $K_5$ come from distinct tripartites.
\item There is one tripartite giving a triangle in the $K_5$ but all the other 7 edges come from distinct tripartites.
\item There are two tripartites that each give a triangle in the $K_5$, these two triangles sharing a single vertex, and all the other 4 edges come from distinct tripartites.
\end{enumerate}
We now delete at least one edge from each of these remaining $K_5$s. This will be done probabilistically and the precise method will be explained later. The deletions in this second round we call \emph{Type 2} deletions. Once they have been performed, the resulting graph is our final graph $G$.

The graph $G$ is $K_5$-free, by definition, but we now have to show that we have not inadvertently destroyed all the triangles in some set of $n^{6/13}(\log n)^c$ vertices. We begin by checking the more basic requirement that the Type 2 deletions destroy only a small proportion of the edges. That is, we check that the expected number of $K_5$s in $G_1$ is less than the expected number of edges (which is already computed to be $n^{21/13}$). To do this, we split into the three cases mentioned above. To calculate the expected number of $K_5$s of type (i), observe that there are at most $n^5$ choices for the vertex set, and $(n^{9/13})^{10}$ choices for the copies of tripartites giving an edge (since there are $n^{9/13}$ tripartites to choose from and we need 10 of them), and the probability that the vertices of the $K_5$ are in these tripartites as prescribed is $(n^{-7/13})^{20}$ (since the probability that a given vertex is in a given tripartite is $n^{-7/13}$), giving that the expected number of these $K_5$ is $n^{15/13}$. Similarly, the expected number of $K_5$s of type (ii) is $n^{5}(n^{9/13})^{8}(n^{-7/13})^{17}=n^{18/13}$. Finally, the expected number of $K_5$s of type (iii) is $n^{5}(n^{9/13})^6(n^{-7/13})^{14}=n^{21/13}$. This last number is roughly equal to the expected number of edges, therefore we will need to modify the parameters by $\log n$ factors. However, the main point is that after this second round of deletions, most edges of the original graph are still present. 

In order to finish off the proof, there are two main difficulties to overcome. The first one is that even though we have made sure that globally not too many edges are deleted, this is, as we have already mentioned, just a necessary condition for the argument to have a chance of working. What we actually need is the stronger statement that every induced subset of size $n^{6/13}(\log n)^c$ still contains a triangle. We can hope that the small set of edges we have removed is ``sufficiently random" for this to be the case, but actually proving that takes some work. Let us sketch how we do it. From now on, it will be convenient to think of each tripartite as having a colour: accordingly, we call the tripartites ``colour classes". If a vertex belongs to, say, the red tripartite, then we say that that vertex is red. 

Let us now fix a set $A$ of size $n^{6/13}(\log n)^c$. As shown above, we can take it for granted that $G_0$ contains a big set $\mathcal{T}$ of triangles in $A$, all coming from different colour classes. Moreover, these triangles will be uniformly distributed over $A$. Let $T_C\in \mathcal{T}$ be a triangle coming from the colour class $C$. (Note that not every colour gives a triangle, and not every triangle in $A$ comes from just one colour class.) Let us first deal with Type 1 deletions. An edge of some $T_C$ gets deleted by the Type 1 deletions if the endpoints of this edge share a colour other than $C$. So intuitively we can imagine that $G_0$ has already been constructed, and then we place these triangles $T_C$ randomly inside $A$ and hope that most triangles will not have any edge contained in another colour class. It is not too hard to show, under suitable assumptions, that with very high probability the density of pairs of vertices in $A$ sharing a colour is fairly low (this essentially comes from the fact that the typical sizes of the tripartites are smaller - after adjusting the parameters by suitable log factors - than the size of $A$). Therefore for a fixed $T_C$ it is indeed true that with fairly high probability its edges will not be deleted by Type 1 deletions. However, these events are not independent for different colours $C$. To overcome this difficulty, we define a set $\Pi$ of roughly $\log n$ partitions with the property that for any pair of distinct colours $C,D$ there is a $\pi\in \Pi$ such that $D$ is in the first part of $\pi$ and $C$ is in the second part. We now define a $\pi$-\emph{dangerous pair} to be a pair of vertices that share a colour from the first part of $\pi$. If an edge $xy$ of a $T_C$ gets deleted (by Type 1 deletions) then $x$ and $y$ share a colour $D\neq C$ and there is some $\pi\in \Pi$ such that $D$ is in the first part of $\pi$ and $C$ is in the second part of $\pi$ and therefore $(x,y)$ is a $\pi$-dangerous pair. But note that, as indicated above, the density of $\pi$-dangerous pairs will be fairly low, so the probability that an edge of $T_C$ is deleted because of a colour in the first part of $\pi$ is low, and, conditional on the outcome of colours in the first part of $\pi$, these events are now independent for all $C$ in the second part of $\pi$. We can therefore conclude that only a small proportion of these $T_C$s will lose an edge thanks to colours in the first part of $\pi$. Thus, since $\Pi$ is small, we deduce that most triangles $T_C$ will not lose an edge. That is, we can find many triangles in $A$ even after the Type 1 deletions.

Now let us define Type 2 deletions. Given the graph $G_1$, we order its edges randomly and keep each edge provided that it does not form a $K_5$ when combined with the edges that we have already decided to keep. We remark that this construction is a variant of the so called $K_5$-free process. The edges we keep will form our final graph $G$.

To be more precise, we note here that in fact we keep an edge only if it does not form a so called \emph{core} of a $K_5$ of $G_1$ when combined with the edges that we have already decided to keep. The core is a certain subgraph of a $K_5$ defined in terms of the colours of its edges. The reader is encouraged to think of the core of a $K_5$ as the $K_5$ itself (especially as we can prove that the core of any $K_t$ is itself, but the proof of this fact is very long and we do not include it in this paper).

As shown above, the number of $K_5$s in $G_1$ is less than the number of edges, that is, on average an edge is contained in less than one $K_5$. In fact, one can show that almost surely \emph{every} edge will be contained in a relatively small number of $K_5$s. It is not hard to see that this means that any triangle in $G_1$ is also present in $G$ with probability not very close to 0. Since the number of triangles in $G_1\lbrack A\rbrack$ is large, standard concentration inequalities will imply that with very high probability $G\lbrack A\rbrack$ still contains a triangle. Using the union bound over all $A$ (of size roughly $n^{6/13}$), we conclude that almost surely every $G\lbrack A\rbrack$ contains a triangle, finishing the proof.

Let us briefly discuss how we determined the parameters of our construction. Let $n^{\delta}$ be the number of tripartite copies placed, let $n^{\beta}$ be the size of each part of each of these copies, and let $n^{\alpha}$ be the set size that will guarantee an induced triangle. The parameters $\delta,\beta$ have been chosen to optimize the result: that is, to allow $\alpha$ to be as small as possible. There are three main conditions that we need to impose on these parameters.

The first one is that we need enough triangles in $G_0$ inside every $A$ of size $n^{\alpha}$. It is not hard to see that this condition is equivalent to
\begin{equation} \label{enoughtriangles}
	\delta+3(\alpha+\beta-1)\geq \alpha.
\end{equation}
The second one comes from the fact that the parts of the tripartites will not contain a triangle in $G$ (since every edge inside a part of a tripartite gets deleted by Type 1 deletions), so we trivially need
\begin{equation} \label{partsnotriangles}
	\alpha\geq \beta.
\end{equation}
Finally, we want the expected number of $K_5$s in $G_1$ to be less than the expected number of edges in $G_1$ which gives (only considering those $K_5$s which are type (iii) in the sense described a few paragraphs above)

\begin{equation} \label{k5freqs}
	\delta+2\beta \geq 5+6\delta+14(\beta-1).
\end{equation}

It is not hard to see that these conditions force $\alpha\geq 6/13$ and that equality is achieved by taking $\delta=9/13,\beta=6/13$.

This leads us to the other main difficulty, which arises only when we consider more general values of $s,t$. While (\ref{enoughtriangles}) is essentially the same but with 3 replaced by $s$, and (\ref{partsnotriangles}) is exactly the same, (\ref{k5freqs}) becomes completely different. Indeed, it will be crucial to analyse all possible ways that a $K_t$ can occur in $G_1$ in some systematic way, rather than writing down the three possibilities (i),(ii),(iii) as we did above in the $s=3,t=5$ case, since in general there are many ways that a $K_t$ can be formed from the contributions of the various $s$-partite graphs. Analysing these decompositions of $K_t$, which we shall refer to as colour schemes (again by imagining that each $s$-partite graph has its own colour), is necessary to determine the best parameters $\delta,\beta$, and also to prove Theorem \ref{ERbound} for these parameters. The complicated formula for $\alpha$ is obtained by solving the system of inequalities (\ref{enoughtriangles}),(\ref{partsnotriangles}),(\ref{k5freqs}) that we obtain in the general case.

\vspace{2mm}

The organization of this paper is as follows. In Section \ref{sectionconstruction} we present our construction. In Section \ref{sectionmain} we give the main part of the proof conditional on three lemmas. These lemmas are proved in Section \ref{sectionauxiliary}. The first one, which asserts that each edge in $G_1$ is contained in a small number of (cores of) $K_t$s, is proved in Subsection \ref{subsectionfewcores}, conditional on a lemma about colour schemes that is proved in Subsection \ref{subsectionnegscheme}. The result that says that $G_1\lbrack A\rbrack$ contains many $K_s$s is proved in Subsection \ref{subsectionfindlarget}. Finally, there is an appendix that contains some tedious computations and the source code of a program relevant to some results in Subsection \ref{subsectionnegscheme}.

\section{The precise construction and the main result} \label{sectionconstruction}

\begin{remark}
	Logarithms throughout the paper are to base $e$. We will not be concerned with floor signs, divisibility, and so on. Also, we will tacitly assume that $n$ is sufficiently large whenever this is needed. Moreover, throughout the rest of the paper, it is to be understood that $s\geq 3$ and that $s+2\leq t\leq 2s-1$. Recall that a pair $(s,t)$ is regular if $s\geq 11$ and $s+3\leq t\leq 2s-4$ or if $(s,t)\in\{(10,14),(10,15)\}$, and otherwise it is exceptional.
\end{remark}

\vspace{2mm}

Let \begin{equation*}
\delta=s-(2s-1)\alpha=\begin{cases}
\delta(1)=\frac{(2s-2)(t-s)(t+s-1)+2s^2-4st+2t+2s}{(2s-3)(t-s)(t+s-1)-2s+4}, & \text{ if } (s,t) \text{ is regular} \\
\delta(2)=\frac{(2s-2)(t-s)(s-1)-st+3s-1}{(2s-3)(t-s)(s-1)+2s-t}, & \text{ if } (s,t) \text{ is exceptional}
\end{cases}
\end{equation*}

\begin{lemma} \label{alphaanddelta}
	$\delta<2\alpha<1$.
\end{lemma}

\vspace{-3mm}

\begin{proof2}
	If $(s,t)$ is regular, then
	\begin{align*} 2\alpha-\delta&=\frac{4st-2s^2+2t-6s-2(t-s)(t+s+1)}{(2s-3)(t-s)(t+s-1)-2s+4} \\
	&=\frac{4st-2t^2-4s}{(2s-3)(t-s)(t+s-1)-2s+4}> 0,
	\end{align*} since $s+1\leq t\leq 2s-2$.
	If $(s,t)$ is exceptional, then
	$$2\alpha-\delta=\frac{st-s-1-2(t-s)(s-1)}{(2s-3)(t-s)(s-1)+2s-t}=\frac{2s^2-st+2t-3s-1}{(2s-3)(t-s)(s-1)+2s-t}>0,$$ since $s+1\leq t\leq 2s-1$.
	
	By Lemma \ref{app2} (e) from the appendix, we have $\delta>2/3>1/2$, which implies that $\alpha<1/2$.
\end{proof2}

\begin{remark}
	 Intuitively, one can think of $\alpha$ as $1/2-\epsilon$ for $\epsilon$ quite small and $\delta=1/2+(2s-1)\epsilon$. This makes $\delta$ significantly greater than 1/2 but less than 1. Also, it may be helpful to bear in mind the case $s=3,t=5$, where, as we have seen, $\delta=9/13$ and $\alpha=6/13$.
\end{remark}

\vspace{3mm}

Let \begin{equation*}
	m=n^{\delta}(\log n)^{-c_1}
\end{equation*}
\begin{equation*}
	\gamma=n^{\alpha-1}(\log n)^{-c_2}
\end{equation*}
\begin{equation*}
	a=n^{\alpha}(\log n)^{c_3}
\end{equation*}
where $c_1,c_2,c_3$ are positive constants, to be specified, that depend on $s$ and $t$. (In fact, $c_1$ can be taken to be 0. All we need are that $c_2$ is suitably large and that $c_3$ is sufficiently larger than $c_1,c_2$.)

The following estimates will be used several times later in the paper.

\begin{lemma} \label{mandgamma}
	$m\gamma>1$ and $m\gamma^2<1$.
\end{lemma}

\begin{proof2}
	Note that $\delta+(\alpha-1)=(s-1)-(2s-2)\alpha>0$ since $\alpha<1/2$. This implies that $m\gamma>1$.
	
	Also, $\delta+2(\alpha-1)<4\alpha-2<0$, by Lemma \ref{alphaanddelta}. This implies that $m\gamma^2<1$.
\end{proof2}

\medskip

We construct the graph $G_0$ as follows. Let $V=V(G_0)=\{1,2,...,n\}$. Define independent random subsets $S_1,...,S_{m}$ of $V$ in such a way that each $S_i$ contains each $v\in V$ independently with probability $\gamma$. We call $S_i$ the $i$th \emph{colour class}. If $v\in S_i$, we say that $v$ \emph{has colour} $i$. Now randomly partition each $S_i$ into $s$ sets, $S_{i1},S_{i2},...,S_{is}$ by placing each element of $S_i$ independently at random in one of these parts, and use these sets to define a complete $s$-partite graph. Let $G_0$ be the union of these $s$-partite graphs. We say that a pair of vertices has colour $i$ if both its members have colour $i$. We do not require the pair to form an edge in $G_0$. Remove all edges of $G_0$ that have at least two colours to obtain the subgraph $G_1$. Again, we do not require both colours to give an edge. Another way to state the condition is that if $xy$ is an edge of colour $i$ and $x$ and $y$ both have colour $j$ for some $j\ne i$, then we remove the edge $xy$ even if $x$ and $y$ belong to the same set $S_{jr}$. Finally, for every $K_t$ in $G_1$ we randomly remove a certain edge, which we shall specify in a moment. The resulting graph is called~$G$. 

The graph $G$ is obviously $K_t$-free. We shall prove that for suitable choices for the constants $c_1,c_2,c_3$, we have the following result, which is our main theorem.

\begin{theorem} \label{mainresult}
	For $n$ sufficiently large, there is a positive probability that every subset $A$ of $G$ with $|A|=a$ contains a $K_s$.
\end{theorem}

\vspace{2mm}

\noindent Obviously Theorem \ref{mainresult} implies Theorem \ref{ERbound}.

\vspace{2mm}

Let us now specify which edges are removed from $G_1$. Suppose that $x_1,...,x_t$ form a $K_t$ in $G_1$. Then necessarily any two distinct vertices $x_i$ and $x_j$ share precisely one colour. Indeed, they must share at least one colour since $x_ix_j\in E(G_0)$ but they cannot share more than one since then $x_ix_j$ would have been removed from $G_0$ during the first round of deletions.

\begin{definition} \label{scheme}
	A \emph{colour scheme} for $K_t$ with parameter $s$, or \emph{scheme} for short, is a set $X$ of $t$ nodes and a set $\mathcal{D}$ of subsets of $X$, which we call \emph{colours}, or \emph{blocks}, such that
	\begin{enumerate}[label=(\roman*)]
	\item For any $x,y\in X$, there is a unique $D\in \mathcal{D}$ such that $x,y$ both belong to $D$.
	\item Every colour appears on at least two nodes.
	\item Every colour appears at most $s$ times.
	\end{enumerate}
	\vspace{2mm}
	
\noindent A pair of nodes is called an \emph{edge} and the \emph{colour} of an edge is the unique colour that contains both endpoints. (Note that a node may have several colours.) If a node $x$ belongs to a colour $D$, we shall say that $D$ \emph{labels} $x$. We also define a \emph{label} to be a pair $(x,D)$ such that $x$ is a node and $D$ labels $x$. The number of labels in a scheme is thus the sum of the sizes of all the colours.
\end{definition}

\vspace{2mm}

If $X=\{x_1,...,x_t\}$ forms a $K_t$ in $G_1$, then there is set of (at most ${t \choose 2}$) colours such that $X$ is a colour scheme with respect to those colours, and no other colour labels more than one vertex in $X$. Indeed, we have already observed that property (i) holds. Choosing the colours suitably, (ii) can clearly be achieved. For property (iii), observe that if some colour $D$ labels at least $s+1$ vertices, then there must exist distinct vertices $x_i$ and $x_j$ that belong to the same part of the complete $s$-partite graph of colour $D$. Then $D$ does not provide an edge between $x_i$ and $x_j$, so some other colour must, but then $x_i$ and $x_j$ share at least two colours, which contradicts (i).

Thus, any $K_t$ in $G_1$ can be viewed as a scheme in a natural way. A simple upper bound for the expected number of $K_t$s associated with a scheme $Q$ is $n^tm^b\gamma^l$, where $l$ is the number of labels of $Q$ and $b$ is the number of colours of $Q$. Indeed, the number of ways choosing the $t$ nodes is at most $n^t$, the number of ways of choosing the $b$ colours (from the $m$ colours used to construct $G_1$) is at most $m^b$, and the probability that any given choice of nodes and colours realizes the scheme is $\gamma^l$, since for each label the probability that the given node receives the given colour is $\gamma$, and all these events are independent.

Now $n^tm^b\gamma^l=n^{t+b\delta+l(\alpha-1)}(\log n)^f$ for some $f=f(s,t,b,l)$. Also, once we know that a certain pair $u,v$ of vertices have a colour in common, the expected number of $K_t$s associated with $Q$ that contain $u$ and $v$ becomes at most roughly $n^{t-2}m^{b-1}\gamma^{l-2}=n^{t-2+(b-1)\delta+(l-2)(\alpha-1)}(\log n)^{f'}$. This motivates the following definition.


\begin{definition}	\label{schemevalue}
	The \textit{value} of a scheme $Q$ with $b$ colours and $l$ labels, denoted $v(Q)$, is given by the formula $$ v(Q)=t-2+(b-1)\delta+(l-2)(\alpha-1). $$
\end{definition}
Thus, roughly speaking, the expected number of $K_t$s associated with a scheme $Q$ that contain a given edge in $G_1$ is at most $n^{v(Q)}$ up to log factors. The following lemma -- proved in Subsection \ref{subsectionnegscheme} -- shows that this number is small.

\begin{lemma} \label{negscheme}
	Let $Q$ be a scheme. Then $v(Q)\leq 0$.
\end{lemma}

We shall also need a generalization of the notion of a scheme where a pair of nodes does not need to have a colour, if it does have a colour then that colour does not have to be unique, and a colour is allowed to label more than $s$ nodes.

\begin{definition} \label{colourconfig}
	A \emph{colour configuration} consists of a set of nodes and a set of colours labelling the nodes such that every colour appears on at least two nodes.
	
	Given a colour configuration $W$ and a subset $S$ of its nodes, we define the \emph{subconfiguration induced by} $S$ to be the configuration whose nodes are the elements of $S$ and whose colours are the colours of $W$ that appear at least twice on $S$ (which then label the nodes in $S$ that they labelled in $W$).
	
	The \emph{value} of a configuration $W$ is defined to be 
	\[v(W)=h-2+(b-1)\delta+(l-2)(\alpha-1),\] 
where $h$ is the number of nodes, $b$ is the number of colours and $l$ is the number of labels in $W$ (where a label is again a pair $(x,D)$ where $x$ is a node labelled by the colour $D$).	
\end{definition}

The same argument as for schemes shows that, once we condition on the event that $u$ and $v$ are both coloured red, the expected number of occurrences of a colour configuration $W$ that contain both $u$ and $v$ is at most $n^{v(W)}$ up to log factors. (In fact, it is smaller unless $u$ and $v$ share a colour in $W$.)


\begin{definition} \label{defncore}
	The \emph{core} of a scheme $Q$, denoted $C(Q)$, is the induced subconfiguration $S$ on at least two nodes for which $v(S)$ is minimal. If several subconfigurations have the same value then the core is the one with the maximum number of nodes. If this is still not unique, then we simply pick an arbitrary one with the given properties.
\end{definition}

\begin{remark}
	We can in fact prove that $C(Q)=Q$ for every scheme $Q$. Although using that fact would simplify the argument in this paper slightly, this gain does not compensate for the extra work needed to establish it, so we shall avoid using it. Nevertheless, the reader is encouraged to think of a core just as a scheme: that is, as a $K_t$ in the graph $G_1$ with the colours given by the $s$-partite graphs with vertex sets that contain at least two of its vertices. 
\end{remark}

\begin{lemma} \label{corefacts}
	Let $Q$ be a scheme. Then $C(Q)$ has at least 3 nodes, $v(C(Q))\leq 0$, and $v(S)\geq v(C(Q))$ for every induced subconfiguration $S$ of $C(Q)$ with at least two nodes.
\end{lemma}

\begin{proof2}
	The first two assertions follow from Lemma \ref{negscheme}, since an induced subconfiguration of $Q$ with two nodes has value 0. The third assertion follows immediately from the definition of the core.
\end{proof2}

\vspace{2mm}

We can now define $G$ precisely. Following an idea in \cite{wolfovitz}, we assign independently to each edge $e$ of $G_1$ a birthtime $\beta_e$, chosen uniformly randomly from $\lbrack 0,1 \rbrack$. Equivalently, we order the edges of $G_1$ uniformly at random from all the possible orderings. To define the edge set $E(G)$, which will be a subset of $E(G_1)$, we recursively decide for each $e\in E(G_1)$ whether $e\in E(G)$, as follows. Suppose that the decision has been made for every $e'\in E(G_1)$ with $\beta_{e'}<\beta_e$. Then let $e\in E(G)$ unless there is a $K_t$ in $G_1$, which we view as a scheme $Q$, for which the edges of $C(Q)$ all have birthtime at most $\beta_e$ and they all (apart from $e$) already belong to $E(G)$.

For any $K_t$ in $G_1$ there is an edge in the core of that $K_t$ that is not an edge of $G$, since if all the edges in the core apart from the last one are chosen to belong to $E(G)$, then the last one is not. Thus, $G$ is $K_t$-free. It remains to prove that with positive probability every set of $a$ vertices still contains a $K_s$, which was Theorem \ref{mainresult} above.

\section{The proof of Theorem \ref{mainresult}} \label{sectionmain}

In this section, we shall prove Theorem \ref{mainresult} conditional on two lemmas, which we shall prove in Section \ref{sectionauxiliary} and which are where most of the work will be. The first one says, roughly speaking, that for any $A$ of size $a$, the induced subgraph $G_1\lbrack A \rbrack$ of $G_1$ contains many copies of $K_s$.
\begin{lemma} \label{findlarget}
	Almost surely, for every $A$ of size $a$ there is a set of $\Omega(ma^s\gamma^s)$ monochromatic copies of $K_s$ inside $G_1\lbrack A\rbrack$, each with a different colour.
\end{lemma}
\noindent The second tells us that any edge in $G_1$ is contained in few cores. Here, and in what follows, we use the word ``core" to refer to the core of a $K_t$ in $G_1$.
\begin{lemma} \label{fewcores}
	Almost surely, any edge in $G_1$ is contained in at most $(\log n)^{2t}$ cores.
\end{lemma}

We shall use McDiarmid's inequality \cite{mcdiarmid} in the next proof, which for convenience we recall here. Let $Y_1,\dots,Y_N$ be independent random variables, taking values in a set $S$, and let $X=g(Y_1,\dots,Y_N)$ for some $g:S^N\rightarrow \mathbb{R}$ with the property that if $y,y'\in S^N$ only differ in their $i$th coordinate, then $|g(y)-g(y')|\leq c_i$. Then the inequality states that 
\[\mathbb{P}\big\lbrack |X-\mathbb{E}\lbrack X\rbrack|\geq r \big\rbrack\leq 2 \exp\Big(\frac{-2r^2}{\sum_i c_i^2}\Big).\]

The following lemma, together with Lemmas \ref{findlarget} and \ref{fewcores} and a union bound, implies Theorem \ref{mainresult}.

\begin{lemma}
	Suppose that $G_1$ is such that any edge in $G_1$ is contained in at most $(\log n)^{2t}$ cores. Let $A$ be a set of vertices of size $a$ such that the induced subgraph $G_1\lbrack A\rbrack$ contains $\Omega(ma^s\gamma^s)$ monochromatic copies of $K_s$, each with a different colour. Then the probability, conditional on the graph $G_1$, that $G\lbrack A\rbrack$ does not contain any $K_s$ is $o\big(\frac{1}{{n \choose a}}\big)$.
\end{lemma}

\begin{proof2}
	Choose $\Omega(ma^s\gamma^s)$ monochromatic copies of $K_s$ in $G_1\lbrack A\rbrack$, all of distinct colours. Let the set of these copies be $\mathcal{T}$. Then by the definition of the first deletion process, the elements of $\mathcal{T}$ are edge disjoint. Let $T\in \mathcal{T}$. Let $E_T$ be the set of all edges of cores that have at least one edge that belongs to $T$, together with the edges of $T$ itself. Clearly, $|E_T|\leq \binom s2+\binom s2(\log n)^{2t}{t \choose 2}\leq (\log n)^{3t}$. Let $B_T$ be the event that the birthtimes of the edges of $T$ precede the birthtimes of all other edges in $E_T$. If $B_T$ occurs, then the only way an edge of $T$ could be deleted from $G_1$ and therefore fail to be present in $G$ is if $T$ itself contains a core of some $K_t$. But note that there is no colour that labels every vertex in a core $C$. Indeed, if there is such a colour, then since all edges in a core belong to $G_1$, there is no other colour appearing at least twice on the node set of $C$, therefore $C$, considered as a colour configuration, has value $h-2+(h-2)(\alpha-1)=(h-2)\alpha$ (where $h$ is the number of nodes in $C$), which contradicts Lemma \ref{corefacts}. It follows that if $B_T$ occurs, then every edge of $T$ is present in $G$. 
	
	For a fixed $G_1$, let $X$ be the number of events $B_T$ that occur over all $T\in \mathcal{T}$. Then $X$ is a random variable with the property that if $X\ne 0$, then there is some $T\in\mathcal{T}$ that belongs to $G\lbrack A\rbrack$. It therefore suffices to prove that $\mathbb{P}\lbrack X=0\rbrack=o\big(\frac{1}{{n \choose a}}\big)$. 
	
	To do this, we apply McDiarmid's inequality when $Y_i$ is the birthtime of the $i$th edge. Since the $T\in\mathcal{T}$ are edge disjoint, and any edge $e$ in $G_1$ is contained in at most $(\log n)^{2t}$ cores, it follows that $e$ is contained in at most $1+(\log n)^{2t}{t \choose 2}\leq (\log n)^{3t}$ of the graphs $E_T$. Hence, changing the birthtime $\beta_e$ of $e$ influences at most $(\log n)^{3t}$ of the events $B_T$. Also, if $e\not \in \cup_{T\in \mathcal{T}}E_T$, then $\beta_e$ does not influence \emph{any} event $B_T$. Thus, by McDiarmid's inequality (with some $N\leq |\mathcal{T}|(\log n)^{3t}$), we get $$ \mathbb{P}\lbrack X=0\rbrack\leq 2 \exp\Big(\frac{-2(\mathbb{E}\lbrack X\rbrack)^2}{|\mathcal{T}|(\log n)^{3t}((\log n)^{3t})^2}\Big). $$
	
	Now note that $\mathbb{P}\lbrack B_T\rbrack\geq |E_T|^{-{s\choose 2}}\geq (\log n)^{-3s^2t}$, so $\mathbb{E}\lbrack X\rbrack \geq |\mathcal{T}|(\log n)^{-3s^2t}$, and 
	\[\mathbb{P}\lbrack X=0\rbrack\leq 2 \exp\Big(\frac{-2|\mathcal{T}|}{(\log n)^{6s^2t+9t}}\Big).\]
	
	Finally, note that ${n \choose a}\leq n^a= \exp(a\log n)$. To finish the proof we just need to verify that $\frac{|\mathcal{T}|}{(\log n)^{6s^2t+9t}}=\omega(a\log n)$. Since $$\frac{|\mathcal{T}|}{a}=\Omega(ma^{s-1}\gamma^s)=n^{\delta+(s-1)\alpha+s(\alpha-1)}(\log n)^{-c_1+(s-1)c_3-sc_2}=(\log n)^{-c_1+(s-1)c_3-sc_2},$$ we are done provided that $(s-1)c_3-sc_2-c_1>6s^2t+9t+1$.
\end{proof2}

\section{The proofs of the auxiliary lemmas} \label{sectionauxiliary}

In this section we shall prove Lemmas \ref{negscheme}, \ref{findlarget} and \ref{fewcores}, which are the results we used in the proof of Theorem \ref{mainresult} but have not yet proved.

\subsection{The proof of Lemma \ref{fewcores}} \label{subsectionfewcores}

Let $e$ be an edge in $G_1$. We would like to show that it belongs to at most $(\log n)^{2t}$ cores. Any core that contains $e$ can be viewed as a core in a scheme that contains $e$, and as such it has nonpositive value. But for any colour configuration $W$ (with more than two labels), the expected number of occurrences of that colour configuration in $G_0$ containing a fixed edge in $G_1$ is at most $n^{v(W)}(\log n)^{-c_2}$ (as we remarked slightly less precisely after Definition \ref{colourconfig}), which is at most $(\log n)^{-c_2}$ if $v(W)\leq 0$. In particular, the probability that an edge $e$ in $G_1$ is contained in $r$ cores that are pairwise disjoint apart from their intersection on $e$ is at most $(\log n)^{-rc_2}$. If $r=\log n$ then this is much less than $1/n^2$, and therefore almost surely no edge is contained in $\log n$ cores of the above form.

In general, the cores containing $e$ need not be disjoint. This adds a complication, and we need to introduce a few definitions to handle it, but the main reason Lemma \ref{fewcores} holds is the one given in the previous paragraph. The next definition describes the kind of colour configuration which -- if it occurs in $G_0$ -- can produce many cores in $G_1$ (that is, cores of $K_t$s in $G_1$ that we view as schemes) that contain a given edge $xy$. Soon we shall argue that almost surely no such large configuration occurs in $G_0$.

\begin{definition} \label{corecontainer}
	An \emph{abstract core container} $W$ is a colour configuration whose nodes are $\{x\}\cup\{y\}\cup Z$ and in which every $z\in Z$ is contained in at least one abstract core, where an abstract core is defined as follows.
	
	An \emph{abstract core} in a core container is an induced subconfiguration $S$ consisting of at most $t$ nodes and containing $x$ and $y$ such that for any induced subconfiguration $S'\subset S$ containing $x,y$, we have $v(S')\geq v(S)$ and such that for any two distinct $u,v\in S$ there is a unique colour that labels both $u$ and $v$.
	
	The \emph{size} of a core container is the number of nodes it contains.
	
	A core container is \emph{irreducible} if it is not possible to remove a label or colour and still have a core container.
\end{definition}

\begin{remark}
	Assume for a moment that we know that the core of a scheme is the scheme itself (see the remark after Definition \ref{defncore}). Then Lemma \ref{fewcores} just asserts that each edge in $G_1$ is contained in few $K_t$s. Then we can replace the technical notion of abstract core container with the notion of \emph{abstract scheme container} instead. What we mean by that is a colour configuration whose nodes are $\{x\}\cup \{y\}\cup Z$ and in which every $z\in Z$ is contained in at least one colour scheme containing $x$ and $y$ as well. This is a configuration that is dangerous to us since if it occurs in $G_0$, then the edge $xy$ is contained in many $K_t$s (corresponding to the various schemes in the configuration).
\end{remark}

\medskip

Note that as the vertices of $G_0$ are coloured, we can naturally talk about $G_0$ containing various colour configurations. We shall now establish that:

\begin{enumerate}
	\item If an edge in $G_1$ is contained in many cores, then there is a large irreducible core container in $G_0$.
	
	\item There are not too many irreducible abstract core containers of fixed size.
	
	\item The expected number of occurrences in $G_0$ of any large abstract core container is small.
\end{enumerate}

The last two points will imply that almost surely there is no large irreducible core container in $G_0$, which in turn implies that there is no edge in $G_1$ that is contained in many cores.

Note that for the second point it is important that we count only irreducible core containers because otherwise the number of colours in the core container could be arbitrarily large.

\begin{lemma} \label{coretocont}
	If the edge $e=uv$ is contained in at least $r$ cores of $K_t$s in $G_1$, then there is an irreducible core container $W$ in $G_0$ with $x=u,y=v$ (as in Definition \ref{corecontainer}) and with size between $\frac{1}{2}r^{1/t}$ and $tr$.
\end{lemma}

\begin{proof2}
	Define a colour configuration $W_0$ as follows. Arbitrarily pick $r$ cores that contain $e$. The set of nodes of $W_0$ is the set of vertices of $G_1$ that are in one of these $r$ cores. The set of colours is the set of those colours in $G_0$ that appear at least twice on this set of nodes. This does indeed define a core container, since any core of a $K_t$ in $G_1$ that contains $e$ satisfies the two properties required of an abstract core in $W_0$: the minimality of $v$ follows from the definition of a core, and the condition about the colours follows from the fact that the $K_t$ belongs to $G_1$.
	
	How many nodes does $W_0$ have? Any core consists of between 2 and $t$ nodes, so if the number of nodes of $W_0$ is $h$, then $r\leq \sum_{2\leq j\leq t} {h \choose j}\leq (2h)^t$. Thus, $h\geq \frac{1}{2}r^{1/t}$. On the other hand, $h\leq rt$, since the vertex set of $W_0$ is a union of $r$ cores. Now remove labels or colours as long as we still get a core container; the object we end up with is an irreducible core container of the required size. 
\end{proof2}

\begin{lemma} \label{containercount}
	The number of distinct irreducible abstract core containers of size $h$ is at most $ht^2\cdot 2^{ht^2}\cdot h^{ht^2}$.	
\end{lemma}

\begin{proof2}
	First we shall prove that the number of labels in an irreducible core container of size $h$ is at most $2h{t \choose 2}\leq ht^2$. For any occurrence of a colour $D$ at some node $u$ (that is, for any label $(u,D)$), there must exist $v\in \{x\}\cup\{y\}\cup Z$ such that every abstract core containing $v$ contains $u$ and the colour $D$, or else we could remove the occurrence of $D$ at $u$ and still have a core configuration. But for any $v$, there are at most $2{t \choose 2}$ such pairs $(u,D)$, since $u$ must belong to the intersection of the vertex sets of the abstract cores containing $v$, and in a given abstract core there are at most $2{t \choose 2}$ labels. Indeed, an abstract core is an induced subconfiguration so each of its colours labels at least two nodes. Now if an abstract core has $q$ colours and they label $d_1,\dots,d_q$ nodes, then $\sum_{i\leq q} {d_i \choose 2}\leq {t \choose 2}$ because the abstract core has at most ${t \choose 2}$ pairs of nodes. Since $d_i\geq 2$ for each $i$, it follows that $\sum_{i\leq q} d_i\leq 2{t\choose 2}$.
		
	So there are at most $ht^2$ choices for the total number of labels. Since the partition function $p(k)$ is at most $2^k$, it follows that for each possibility for the number of labels, there are at most $2^{ht^2}$ choices for the number of occurrences for each colour class. Suppose we have $b$ colours and the numbers of times that they occur are $l_1,\dots,l_b$. Then the number of choices for the vertices labelled by these colours is at most ${h \choose l_1}{h \choose l_2}\dots {h \choose l_b}\leq h^{l_1+\dots+l_b}\leq h^{ht^2}$.	
\end{proof2}

\medskip

Next, we shall investigate how many copies we expect to have in $G_0$ of a given abstract core container. Let $W$ be more generally any colour configuration with $h$ nodes, $b$ colours and $l$ labels. Then the expected number of occurrences of such a configuration is at most $n^hm^b\gamma^l$. Indeed, the number of ways of choosing the $h$ nodes is at most $n^h$. The number of ways of choosing the $b$ colours is at most $m^b$. And for each label, the probability that the given node receives the given colour is $\gamma$, and all these events are independent, so the probability that any given choice of nodes and colours realizes the scheme is $\gamma^l$.

\begin{definition}
	We call $n^hm^b\gamma^l$ the \emph{frequency} of the configuration $W$ and denote it by $\omega(W)$.
\end{definition}

\begin{lemma} \label{contfreq}
	Let $W$ be an abstract core container of size $h$. Then 
	\[\omega(W)\leq n^2(\log n)^{-\frac{h-2}{t}c_2}.\]
\end{lemma}

To prove this result, we will kill some of the nodes and colours and remove some of the labels of the core container in steps. To keep track of which nodes and colours have been killed, we introduce the following definition.

\begin{definition} \label{partialconfig}
	A \emph{partial configuration} $P$ consists of four pairwise disjoint sets $\{x\}$, $\{y\}$, $Z_0$ and $Z_1$ of nodes, and two disjoint sets $\mathcal{B}_0,\mathcal{B}_1$ of colours that label those nodes in such a way that any $B\in \mathcal{B}_1$ labels at least two nodes. We write $\mathcal{B}$ for $\mathcal{B}_0\cup \mathcal{B}_1$ and $Z$ for $Z_0\cup Z_1$.
	\end{definition}
	
	We now generalize the notion of frequency to this setting, which can be thought of as the expected number of occurrences of the colour configuration for given choices of the nodes in $Z_0$ and colours in $\mathcal{B}_0$, which represent the nodes and colours that have already been killed. Thus, we let $r=|\{x\}\cup \{y\}\cup Z_1|$ be the number of nodes yet to choose, we let $g=|\mathcal{B}_1|$ be the number of colours yet to choose, and we let $u$ be the total number of labels, including the labels on nodes in $Z_0$ and of colours in $\mathcal{B}_0$. Then we can choose the remaining nodes in at most $n^r$ ways and the remaining colours in at most $m^g$ ways, and for each label there is a probability $\gamma$ that the given node receives the given colour. So we define the frequency $\omega(P)$ to be $n^rm^g\gamma^u$.

\begin{proof}[\textnormal{\textbf{Proof of Lemma \ref{contfreq}.}}]
	We shall define a sequence $P_0,\dots,P_k$ of partial configurations such that $\omega(P_0)=\omega(W)$, $\omega(P_k)\leq n^2$, $k\geq \frac{h-2}{t}$ and $\omega(P_j)\geq \omega(P_{j-1})(\log n)^{c_2}$. Clearly, this suffices to prove the lemma.
	
	We shall define the $P_j$ recursively. In what follows we use the notation of Definition  \ref{corecontainer} and Definition \ref{partialconfig}. When there is ambiguity, we will write $Z_0(P)$ to mean $Z_0$ in the partial configuration $P$, and similarly for $Z_1,\mathcal{B}_0,\mathcal{B}_1$. The set of all nodes (respectively, colours) for every $P_j$ will be the same as the set of all nodes (respectively, colours) of $W$, namely $\{x\}\cup \{y\}\cup Z$ (respectively, $\mathcal{B}$). However, $\mathcal{B}_0,\mathcal{B}_1,Z_0,Z_1$ and the labels will be different for the various $P_j$.
	
	Let us define $P_0$ to be the partial configuration whose nodes, colours and labels are the same as those of $W$ and which has $Z_0=\mathcal{B}_0=\emptyset$. Then $\omega(P_0)=\omega(W)$.
	
	Given $P_{j-1}$ with $Z_1(P_{j-1})\neq \emptyset$, we define $P_j$ as follows. Pick some $z\in Z_1(P_{j-1})$ arbitrarily. As $W$ is a core container, we can choose an abstract core $S$ in $W$ that contains $z$. Let $S_1=S\cap Z_1(P_{j-1})$. Let $\mathcal{D}$ be the set of those colours $B\in \mathcal{B}_1(P_{j-1})$ that occur at least twice on $S$ in $P_{j-1}$. Then let the sets of nodes of $P_j$ be $Z_0(P_j)=Z_0(P_{j-1})\cup S_1$ and $Z_1(P_j)=Z_1(P_{j-1})\setminus S_1$, and let the sets of colours be $\mathcal{B}_0(P_j)=\mathcal{B}_0(P_{j-1})\cup \mathcal{D}$ and $\mathcal{B}_1(P_j)=\mathcal{B}_1(P_{j-1})\setminus\mathcal{D}$. The labels of $P_j$ are those of $P_{j-1}$ except that all occurrences of colours in $\mathcal{B}_0(P_j)$ are removed from $S$. It is clear that $P_j$ is a partial configuration.
	
	We want to prove that $\omega(P_j)\geq \omega(P_{j-1})(\log n)^{c_2}$.
	
	\noindent \emph{Claim.} $\frac{\omega(P_j)}{\omega(P_{j-1})}\geq \frac{\omega(S\setminus S_1)}{\omega(S)}$, where $S$ and $S\setminus S_1$ are identified with their induced subconfigurations from $W$.
	
	\medskip
	
	\noindent \emph{Proof of Claim.} The contribution of the nodes is (a factor of) $n^{-|S_1|}$ to both $\frac{\omega(P_j)}{\omega(P_{j-1})}$ and $\frac{\omega(S\setminus S_1)}{\omega(S)}$. Hence it suffices to prove that the contribution of any colour (and its labels) to $\frac{\omega(P_j)}{\omega(P_{j-1})}$ is at least as much as its contribution to $\frac{\omega(S\setminus S_1)}{\omega(S)}$. There are two cases to consider.
	
	\emph{Case 1.} If $B$ is a colour that occurs at most once on $S$ in $W$, then its contribution to $\frac{\omega(S\setminus S_1)}{\omega(S)}$ is 1, whereas its contribution to $\frac{\omega(P_j)}{\omega(P_{j-1})}$ is at least 1. (Indeed, since $m\gamma^2<1$, the contribution of \emph{any} colour to $\frac{\omega(P_j)}{\omega(P_{j-1})}$ is at least 1.)
	
	\emph{Case 2.} Suppose, then, that $B$ is a colour that occurs at least twice on $S$ in $W$. 
	
	\emph{Case 2a.} If $B\in \mathcal{B}_0(P_{j-1})$, then let $d$ be the number of occurrences of $B$ on $S_1$ in $W$. The contribution of $B$ to $\frac{\omega(S \setminus S_1)}{\omega(S)}$ is at most $\gamma^{-d}$. Indeed, this is clear unless $B$ occurs exactly once on $S\setminus S_1$ in $W$. But if this is the case, then the contribution of $B$ is precisely $m^{-1}\gamma^{-(d+1)}$, which is at most $\gamma^{-d}$, by Lemma \ref{mandgamma}.
	
	Note that any node in $S_1$ (and in fact more generally in $Z_1(P_{j-1})$) that is labelled by $B$ in $W$ is also labelled by $B$ in $P_{j-1}$. Therefore, the contribution of $B$ to $\frac{\omega(P_j)}{\omega(P_{j-1})}$ is at least $\gamma^{-d}$.
	
	\emph{Case 2b.} If $B\in \mathcal{B}_1(P_{j-1})$, then let $d$ be the number of occurrences of $B$ on $S$ in $W$. The contribution of $B$ to $\frac{\omega(S \setminus S_1)}{\omega(S)}$ is at most $m^{-1}\gamma^{-d}$. Indeed, this is clear unless $B$ occurs at least twice on $S\setminus S_1$ in $W$. But in this case the contribution of $B$ is at most $\gamma^{-(d-2)}$, which is at most $m^{-1}\gamma^{-d}$, by Lemma \ref{mandgamma}.
	
	Note that any node that is labelled by $B$ in $W$ is also labelled by $B$ in $P_{j-1}$. Therefore, $B\in \mathcal{D}$ and the contribution of $B$ to $\frac{\omega(P_j)}{\omega(P_{j-1})}$ is precisely $m^{-1}\gamma^{-d}$.
	
	This completes the proof of the claim.
	
	\medskip
	
	Since $S$ is an abstract core in $W$, we have $v(S)\leq v(S\setminus S_1)$, by the minimality of $S$. Because $S_1\neq \emptyset$, and every node in a core has a label on it, it follows that, considering $S$ and $S\setminus S_1$ as induced subconfigurations of $W$, we have $\omega(S \setminus S_1)\geq \omega(S)(\log n)^{c_2}$. Using the claim above, the inequality $\omega(P_j)\geq \omega(P_{j-1})(\log n)^{c_2}$ follows.
	
	Eventually we obtain a partial configuration $P_j$ with $Z_1(P_j)=\emptyset$. When this happens, we set $k=j$. By definition, we have in that case that $\omega(P_k)=n^2m^g\gamma^u$ where $g=|\mathcal{B}_1(P_k)|$ and $u$ is the number of labels in $P_k$. Since any $B\in \mathcal{B}_1(P_k)$ labels at least two nodes in $P_k$ and $m\gamma^2\leq 1$, we find that $\omega(P_k)\leq n^2$. Also note that $|Z_1(P_j)|\geq |Z_1(P_{j-1})|-t$ for any $j$, and $|Z_1(P_0)|=|Z|=h-2$, so $k\geq \frac{h-2}{t}$.
\end{proof}

\medskip

We are now in a position to complete the proof of Lemma \ref{fewcores}.

\begin{proof}[\textnormal{\textbf{Proof of Lemma \ref{fewcores}.}}]
	By Lemma \ref{coretocont}, it suffices to prove that in $G_0$ the expected number of irreducible core containers of size between $\log n$ and $(\log n)^{3t}$ is $o(1)$.
	
	\noindent \emph{Claim.} If $\log n \leq h\leq (\log n)^{3t}$, then the expected number of irreducible core containers of size $h$ in $G_0$ is at most $n^2(\log n)^{-t^3h}$.
	
	\noindent \emph{Proof of Claim.} By Lemmas \ref{containercount} and \ref{contfreq}, the expected number of irreducible core containers of size $h$ in $G_0$ is at most $ht^2 2^{ht^2}h^{ht^2}n^2(\log n)^{-\frac{h-2}{t}c_2}\leq h^{3ht^2}n^2(\log n)^{-\frac{h-2}{t}c_2}$. If $c_2\geq 11t^4$, then this is at most $h^{3ht^2}n^2(\log n)^{-11(h-2)t^3}\leq h^{3ht^2}n^2(\log n)^{-10ht^3}\leq n^2(\log n)^{-ht^3}$ so the claim is proved.
	
	\smallskip
	
	But $\sum_{h\geq \log n} n^2(\log n)^{-t^3h}=o(1)$, and the proof is complete.
\end{proof}

\subsection{The proof of Lemma \ref{findlarget}} \label{subsectionfindlarget}

Our proof is based on the following two observations.

\begin{enumerate}
	\item For any set of vertices $A$ of size $a$, $G_0\lbrack A\rbrack$ contains many monochromatic $s$-cliques with pairwise distinct colours.
	
	\item If a monochromatic $s$-clique is present in $G_0$, then it is present also in $G_1$ with high probability, and, crucially, the events that various $s$-cliques are preserved are ``sufficiently independent".
\end{enumerate}

First, we shall construct a small set of bipartitions of the set of colours with a suitable property. In a moment it will become clear why we need this. We will refer to the two parts of a bipartition as the \textit{first part/first half} and the \textit{second part/second half}.

\begin{lemma}
	There exists a constant $c$ and a set $\Pi$ of $c\log n$ partitions of the set of $m$ colours, each into two sets of size $m/2$, such that for any two distinct colours $C$ and $D$ there is a $\pi\in \Pi$ such that $D$ is contained in the first part of $\pi$ and $C$ is contained in the second part of $\pi$.
\end{lemma}

\begin{proof2}
	Take $l=c\log n$ random partitions. For any $C,D$, the probability that none of the partitions is suitable is less than  $(1-\frac{1}{5})^l=n^{-c\log (5/4)}$. For $c$ sufficiently large this is less than $n^{-2}$, which is in turn less than $m^{-2}$ and the result follows from the union bound over all choices of $C,D$.
\end{proof2}

\medskip

Let $xy$ be an edge in $G_0$. Recall that it is not an edge in $G_1$ if $x,y$ have at least two colours in common. Suppose that this is the case. Then there exists some $\pi\in \Pi$ such that $x$ and $y$ have a colour in common from the first half of $\pi$ and also a colour in common from the second half of $\pi$.

\begin{remark}
	From now on, when we say ``the first $m/2$ colours", we will mean ``the $m/2$ colours in the first part of $\pi$" provided it is clear which $\pi$ we are talking about.
\end{remark}

\begin{definition}	
	A pair $(x,y)$ of vertices is \textit{$\pi$-dangerous} for some $\pi\in \Pi$ if there is a colour class among the first $m/2$ colours that contains both $x$ and $y$.
\end{definition}

Fix a set $A$ of vertices with $|A|=a$. Let $\mathcal{D}$ be the collection of colours $D$ such that at least one $K_s$ inside $A$ is entirely coloured with colour $D$ in $G_0$. (We require that every edge is given by this colour: that is, the vertices of the $K_s$ are in different parts of the complete $s$-partite graph with colour $D$.) For each $\pi \in \Pi$, let $\mathcal{D}_{\pi}$ be the set of all $D\in \mathcal{D}$ such that $D$ is one of the last $m/2$ colours. 

To make sense of the statement of the next lemma, the reader should recall that $a\gamma$ is significantly less than 1. (See the beginning of Section \ref{sectionconstruction} for their precise values.)

\begin{lemma} \label{largedpi}
	With probability $1-o(\frac{1}{{n \choose a}})$, $|\mathcal{D}_{\pi}|=\Omega(ma^s\gamma^s)$ for every $\pi \in \Pi$.
\end{lemma}

\begin{proof2}
	Let $C$ be any colour class. The probability that $C$ intersects $A$ in exactly $s$ elements is
\[\mathbb{P}[\Bin(a,\gamma)=s]=
	{a \choose s}\gamma^s(1-\gamma)^{a-s}=\Omega(a^s\gamma^s(1-\gamma)^a)=\Omega(a^s\gamma^s(e^{-2\gamma})^a)=
	\Omega(a^s\gamma^s),\] 
where the last inequality follows from the fact that $a\gamma=n^{2\alpha-1}(\log n)^{c_3-c_2}=o(1)$.
	
	Hence $\mathbb{P}[C\in \mathcal{D}]=\Omega(a^s\gamma^s)$. Moreover, the events $\{C\in \mathcal{D}\}$ are independent. Thus, for any $\pi$, by the Chernoff bound we get $\mathbb{P}\big[|\mathcal{D}_{\pi}|=o(ma^s\gamma^s)\big]\leq e^{-\Omega(ma^s\gamma^s)}$. Therefore, using the union bound over all $\pi\in \Pi$, it suffices to prove that $(\log n)e^{-\Omega(ma^s\gamma^s)}=o(\frac{1}{{n \choose a}})$.
	
	But ${n \choose a}\leq n^a=e^{a\log n}$. Hence, we need $(\log n)e^{-\Omega(ma^s\gamma^s)}=o(e^{-a\log n})$. For this, it is enough to prove that $a\log n=o(ma^s\gamma^s)$, ie. $\log n=o(ma^{s-1}\gamma^s)$. Since
	\begin{equation}
		ma^{s-1}\gamma^s=n^{\delta+(s-1)\alpha+s(\alpha-1)}(\log n)^{-c_1+(s-1)c_3-sc_2}=(\log n)^{-c_1+(s-1)c_3-sc_2}, \label{eqnmagamma}
	\end{equation}
we are done provided that $(s-1)c_3-sc_2-c_1>1$.
\end{proof2}

\vspace{2mm}

Therefore, using the union bound over all sets $A$ of size $a$, we may assume that $|\mathcal{D}_{\pi}|= \Omega(ma^s\gamma^s)$ for every $\pi \in \Pi$ and every such set $A$.

\begin{lemma} \label{dangerousdensity}
	With probability $1-o(1)$ the following holds. For every $A$ of size $a$ and for every $\pi\in \Pi$, the density of $\pi$-dangerous pairs in $A$ is $o(\frac{1}{\log n})$.
\end{lemma}
This result, which we shall prove later, allows us to assume for our fixed set $A$ that the following statement holds.
\begin{enumerate}[label={}]
\item $(\star)$ For any $\pi\in \Pi$, the density of $\pi$-dangerous pairs in $A$ is $o(\frac{1}{\log n})$.
\end{enumerate}

For each $C\in\mathcal D$, pick a $K_s$ uniformly at random in $G_0\lbrack A\rbrack$ of colour $C$, and call it $T_C$. We can now prove that with sufficiently high probability, most $T_C$ will be present in~$G_1$.

\begin{lemma}
	Let $\pi \in \Pi$. Then with probability $1-o(\frac{1}{(\log n){n \choose a}})$, the number of colours $C\in \mathcal{D}_{\pi}$ for which $T_C$ has a $\pi$-dangerous pair of vertices is $o(\frac{|\mathcal{D}_{\pi}|}{\log n})$.
\end{lemma}

\begin{proof2}
	We condition everything on the already chosen first $m/2$ colour classes. Now let $C\in \mathcal{D}_{\pi}$. (Recall that this means that there is a $K_s$ in $A$ in the graph $G_0$ with all its edges of colour $C$, and moreover that $C$ is one of the last $m/2$ colours with respect to $\pi$.) Label the vertices of $T_C$ by $1,2,...,s$. Note that any pair of vertices in $A$ is chosen with equal probability and, by condition $(\star)$, at most $o(\frac{|A|^2}{\log n})$ of them are $\pi$-dangerous. So the probability that the first two vertices of $T_C$ form a $\pi$-dangerous pair is $o(\frac{1}{\log n})$. Hence, for any $C\in \mathcal{D}_{\pi}$, the probability that $T_C$ has a pair of vertices which form a $\pi$-dangerous pair is bounded above by some $p=o(\frac{1}{\log n})$. Moreover, this holds for all such $C$ \textit{independently} of the others. Thus, the probability that $T_C$ contains a $\pi$-dangerous pair for more than $\Omega(\frac{|\mathcal{D}_{\pi}|}{\log n})$ choices of $C\in \mathcal{D}_{\pi}$ is at most $\mathbb{P}\big[\Bin(|\mathcal{D}_{\pi}|,p)=\Omega(\frac{|\mathcal{D}_{\pi}|}{\log n})\big]$. But this is $e^{-\Omega(\frac{|\mathcal{D}_{\pi}|}{\log n})}$. So it remains to show that $(\log n){n \choose a}=o(e^{\Omega(\frac{|\mathcal{D}_{\pi}|}{\log n})})$. Since ${n \choose a}\leq n^a= e^{a\log n}$, it suffices to prove that $a\log n=o(\frac{|\mathcal{D}_{\pi}|}{\log n})$. But $|\mathcal{D}_{\pi}|= \Omega(ma^s\gamma^s)$ so it is enough to prove that $(\log n)^2= o(ma^{s-1}\gamma^s)$. By equation (\ref{eqnmagamma}), this holds provided that $(s-1)c_3-sc_2-c_1>2$.
\end{proof2}

\begin{corollary}
	With probability $1-o(\frac{1}{{n \choose a}})$, for all but $o(|\mathcal{D}|)$ colours $C\in \mathcal{D}$, all the edges of $T_C$ are present in $G_1$.
\end{corollary}

\begin{proof2}
	Suppose that $C\in \mathcal{D}$ and $T_C$ has an edge $e$ which is not present in $G_1$. Then there exists some $\pi\in \Pi$ such that $C$ is in the second half of $\pi$ (so $C\in \mathcal{D}_{\pi}$) and $e$ is $\pi$-dangerous. But by the previous lemma, with probability $1-o(\frac{1}{{n \choose a}})$ the number of such colours $C$ is $o(|\Pi|\cdot \frac{|\mathcal{D}|}{\log n})=o(|\mathcal{D}|)$.
\end{proof2}

\vspace{2mm}

Using Lemma \ref{largedpi} and the union bound over all $A$, Lemma \ref{findlarget} follows.

We now return to proving Lemma \ref{dangerousdensity}. Recall that we want to show that almost surely for every $A$ and every $\pi$, the density of $\pi$-dangerous pairs in $A$ is $o(\frac{1}{\log n})$. This is essentially best possible, since if we choose $A$ to contain one of our colour classes entirely (for a colour chosen from the first part of $\pi$), then the pairs of vertices in that colour class will all be $\pi$-dangerous. Moreover, as the typical size of a colour class is $n\gamma=n^{\alpha}(\log n)^{-c_2}=a(\log n)^{-c_2-c_3}$, the set of these pairs will have density roughly $(\log n)^{-2c_2-2c_3}$.

Accordingly, the next lemma is to make sure that no colour class is exceptionally large.

\begin{lemma}
	With probability $1-o(1)$, the size of every colour class is at most $2n\gamma$.
\end{lemma}

\begin{proof2}
	$\mathbb{P}[\Bin(n,\gamma)> 2n\gamma]=e^{-\Omega(n\gamma)}=o(\frac{1}{m})$. The result follows from the union bound over all colours.
\end{proof2}

\vspace{2mm}

So we may assume that all colour classes have size at most $2n\gamma$.

After applying the union bound over all $\pi\in \Pi$ and $A$, the next result completes the proof of Lemma \ref{findlarget}.

\begin{lemma}
	Fix $\pi\in \Pi$ and a set $A$ of size $a$. With probability $1-o(\frac{1}{(\log n){n \choose a}})$, the number of pairs in $A$ which are $\pi$-dangerous is at most $4\frac{a^2}{(\log n)^{2}}$.
\end{lemma}

\begin{proof2}
	The number of $\pi$-dangerous pairs in $A$ is at most
	\begin{equation} \label{type1bound}
	\sum_{i=1}^{m}\big(\min\{\Bin(a,\gamma),2n\gamma\}\big)^2.
	\end{equation}
Let $h=\frac{a}{m^{1/2}\log n}$. Note that $\log h=(\alpha-\frac{1}{2}\delta)\log n+O(\log \log n)$ and recall that $\alpha>\frac{1}{2}\delta$. Now let $p=\mathbb{P}(\Bin(a,\gamma)\geq h)\leq {a \choose h}\gamma^h\leq (a\gamma)^h\leq e^{-\Omega(h\log n)}$.
	
	Pick some tiny positive $\rho>0$. Note that 
\begin{align*}\mathbb{P}[\Bin(m,p)\geq m^{1/2+\rho}]&\leq {m \choose m^{1/2+\rho}}p^{m^{1/2+\rho}}\leq (mp)^{m^{1/2+\rho}}=e^{-\Omega(m^{1/2+\rho}h\log n)}\\
&=e^{-\Omega(am^{\rho})}= o\bigg(\frac{1}{(\log n){n \choose a}}\bigg).\\
\end{align*}
Therefore we may assume that at most $m^{1/2+\rho}$ of the random variables $\Bin(a,\gamma)$ take value more than $h$.
	
The total contribution to (\ref{type1bound}) of the terms with $\Bin(a,\gamma)\leq h$ is at most $mh^2=\frac{a^2}{(\log n)^{2}}$. The random variable $X\sim \Bin(a,\gamma)$, conditional on $X\geq h$, is bounded above by $h+X'$ where $X'$ is an independent instance of $\Bin(a,\gamma)$. As we assume that all colour classes have size at most $2n\gamma$, it follows that the total contribution to (\ref{type1bound}) of the terms with $\Bin(a,\gamma)\geq h$ is bounded above by
	\begin{equation} \label{type1bound2}
	\sum_{i=1}^{m^{1/2+\rho}} \Big(h+\min\{\Bin(a,\gamma),2n\gamma\}\Big)^2
	\end{equation}
	and we just need to show that this sum is less than $3\frac{a^2}{(\log n)^{2}}$ with probability $1-o(\frac{1}{{m \choose m^{1/2+\rho}}(\log n){n \choose a}})$.
	
	The sum in (\ref{type1bound2}) is at most $m^{1/2+\rho}h^2+(2h+2n\gamma)\sum_{i=1}^{m^{1/2+\rho}}\Bin(a,\gamma)$. The first term is at most $\frac{a^2}{(\log n)^{2}}$. Also, $\log(n\gamma)=\alpha \log n+O(\log \log n)$ and therefore $n\gamma\geq h$, so we just need to show that $\sum_{i=1}^{m^{1/2+\rho}}\Bin(a,\gamma)\leq \frac{a^2}{2n\gamma(\log n)^{2}}$ with the required probability. But the left-hand side is $\Bin(m^{1/2+\rho}a,\gamma)$ and $\mathbb{P}\big[\Bin(m^{1/2+\rho}a,\gamma)\geq \frac{a^2}{2n\gamma(\log n)^{2}}\big]= e^{-\Omega(\frac{a^2}{2n\gamma(\log n)^{2}})}$ since $m^{1/2+\rho}a\gamma = o(\frac{a^2}{2n\gamma(\log n)^{2}})$. This last inequality holds because
	\begin{equation*}
		\log (m^{1/2+\rho}a\gamma)=\big((1/2+\rho)\delta+\alpha+(\alpha-1)\big)\log n+O(\log \log n)
	\end{equation*}
	and	
	\begin{equation*}
		\log(\frac{a^2}{2n\gamma(\log n)^{2}})=\big(2\alpha-1-(\alpha-1)\big)\log n+O(\log\log n),
	\end{equation*}	
	and $(1/2+\rho)\delta+\alpha<1$ for $\rho$ sufficiently small (since $\delta<1$ and $\alpha<1/2$).
	
	Finally, ${m \choose m^{1/2+\rho}}(\log n){n \choose a}= e^{O(a\log n)}$ because $m^{1/2+\rho}=o(a)$ for $\rho$ sufficiently small (as $\delta<2\alpha$). But $a\log n=o(\frac{a^2}{n\gamma(\log n)^{2}})$ provided that $c_3+c_2>3$, so we are done.
\end{proof2}

\subsection{The proof of Lemma \ref{negscheme}} \label{subsectionnegscheme}

It is convenient to introduce the parameter
\begin{equation*}
\eta=2(1-\alpha)-\delta=\begin{cases}
\eta(1)=\frac{-2s^2+4st-2s-6t+8}{(2s-3)(t-s)(t+s-1)-2s+4}, & \text{ if } (s,t) \text{ is regular} \\
\eta(2)=\frac{st-s-2t+3}{(2s-3)(t-s)(s-1)+2s-t}, & \text{ if } (s,t) \text{ is exceptional}
\end{cases}
\end{equation*}

\begin{remark}
	$-\eta$ is the contribution of a block of size two to the value of a scheme. By Lemma \ref{alphaanddelta}, we have $\eta>2-4\alpha>0$.
\end{remark}

\smallskip

The next lemma follows easily from Definition \ref{schemevalue} and is a convenient way to look at the value of a scheme.

\begin{lemma} \label{valueconvenient}
	Let $Q$ be a scheme. Then
	\begin{equation*}
	v(Q)=t+\sum_{D\in \mathcal{D}} (\delta+|D|(\alpha-1))-(\delta+2\alpha)
	\end{equation*}
	where $\mathcal{D}$ is the set of colours in $Q$ and $|D|$ is the number of nodes in $Q$ that are coloured with $D$.
\end{lemma}
	
\vspace{2mm}

We shall now identify a scheme for which equality in Lemma \ref{negscheme} will hold: the value of $\alpha$ was chosen so that the value of this scheme would be 0. This is the (in)equality that generalizes equation (\ref{k5freqs}) from the introduction. This ``extremal scheme" turns out to be different in the regular and the exceptional case, which is why the formula for $\alpha$ also differs in the two cases.

\begin{definition}
	Let $Q_1$ be the scheme where one colour gives a block of size $s$ and the rest of the edges are given by pairwise distinct colours.
	
	Let $Q_2$ be the scheme where one colour gives a block of size $s$, another gives a block of size $t-s+1$ sharing a single vertex with the previous block and the rest of the edges are given by pairwise distinct colours.
\end{definition}

\begin{lemma} \label{extremal}
\begin{enumerate}[label=(\alph*)]
\item If $(s,t)$ is regular, then $v(Q_1)=0$.
\item If $(s,t)$ is exceptional, then $v(Q_2)=0$.
\item If $(s,t)$ is regular, then $v(Q_2)\leq 0$.
\item If $(s,t)$ is exceptional, then $v(Q_1)\leq 0$.
\end{enumerate}
\end{lemma}

\begin{proof2}
We have 
\[v(Q_1)=t+(\delta+s(\alpha-1))+({t \choose 2}-{s \choose 2})(\delta+2(\alpha-1))-(\delta+2\alpha),\] 
and (a) follows by direct substitution.

We also have 
\begin{align*}v(Q_2)=t&+(\delta+s(\alpha-1))+(\delta+(t-s+1)(\alpha-1))\\
&+\Bigl({t \choose 2}-{s \choose 2}-{t-s+1 \choose 2}\Bigr)(\delta+2(\alpha-1))-(\delta+2\alpha),\\
\end{align*}
and (b) follows by direct substitution.
	
	The difference between $Q_1$ and $Q_2$ is that the former contains ${t-s+1 \choose 2}$ edges of distinct colours where the latter contains a block of size $t-s+1$. Using Lemmas \ref{app1} and \ref{app2} (a) from the appendix, we obtain statements (c) and (d). 
\end{proof2}

\begin{definition}
	We call a block in a scheme \emph{large} if it has size at least 3 and \emph{small} otherwise. We call it an $s$-\emph{block} if it has size $s$.
\end{definition}

\vspace{2mm}

We shall begin by proving Lemma \ref{negscheme} in the special case when there is an $s$-block in the scheme.

\begin{lemma} \label{copycase}
	If $Q$ is a scheme and it has an $s$-block then $v(Q)\leq 0$.
\end{lemma}

\begin{proof2}
	Assume that $Q$ is such that $v(Q)$ is maximal. It is enough to show that $Q=Q_1$ or $Q=Q_2$. Since $Q$ has an $s$-block, any other block must have size at most $t-s+1$. By Lemmas \ref{app1} and \ref{app2} (c) from the appendix, any large block of size smaller than $t-s$ gives a smaller contribution to the value than one obtains if the corresponding edges have pairwise distinct colours. Therefore, we may assume that $Q$ has no such block. So every block in $Q$, other than the one of size $s$, has size $2,t-s$ or $t-s+1$. If there is a block of size $t-s+1$, then $Q=Q_2$. If there are no large blocks, then $Q=Q_1$. Otherwise, there is a block of size $t-s\geq 3$.
	
	If there are no other large blocks, then we claim that $v(Q)\leq v(Q_1)$ or $v(Q)\leq v(Q_2)$. Indeed, the $(t-s)$-block can be modified to become a $(t-s+1)$-block (and $Q$ then becomes $Q_2$) and this increases the value provided that $(\alpha-1)\geq (t-s)(-\eta)$, or equivalently $(t-s)\eta\geq (1-\alpha)$. So we may assume that $(t-s)\eta<(1-\alpha)$. But $\delta=s-(2s-1)\alpha>1-\alpha$, since $\alpha<1/2$. Hence, $(t-s)\eta<\delta$, but then $v(Q)\leq v(Q_1)$ by Lemma \ref{app1} from the appendix. 
	
	We may therefore assume that there are at least two large blocks other than the one of size $s$, and that both have size $t-s$. This forces $t-s$ to equal $3$. Moreover, by Lemmas \ref{app1} and \ref{app2} (b), we have that $t=2s-1$. It follows that $s=4$ and $t=7$. So $Q$ consists of a 4-block and several 3-blocks (there can be at most 3) and the rest of the edges are given by distinct colours. It is easy to check that in this case $v(Q)\leq 0$.
\end{proof2}

\vspace{2mm}

Using the previous result, to prove Lemma \ref{negscheme}, it is sufficient to prove the following statement.

\begin{lemma} \label{reduction}
	Suppose that $Q$ is a scheme with $v(Q)$ as large as possible. Assume also that $Q$ does not contain a block of size $s$. Then $v(Q)\leq 0$.	
\end{lemma}

To prove Lemma \ref{reduction}, we shall introduce the following definition.

\begin{definition}
	Let $P$ be a node in a scheme. The \textit{local value} at $P$, which we denote by $v(P)$, is defined by the formula
	\begin{equation*}
	v(P)=1+\sum_{D:P\in D}(\delta/|D|+(\alpha-1)),
	\end{equation*}
where the summation is over all blocks containing $P$.
\end{definition}

\begin{example}
	If $P$ is in a block of size 2 and two blocks of size 4, then 
	\[v(P)=1+3(\alpha-1)+\delta/2+2\cdot\delta/4.\]
\end{example}

\begin{lemma} \label{locsum}
	For any scheme $Q$, we have
	\begin{equation*}
	\sum_P v(P)=v(Q)+(\delta+2\alpha)
	\end{equation*}
	where the summation is over all nodes of $Q$.
\end{lemma}

\begin{proof2}
	This statement follows easily from Lemma \ref{valueconvenient}.
\end{proof2}

\vspace{2mm}

The next result is the key part in the proof of Lemma \ref{negscheme}.

\begin{lemma} \label{localneg}
	Suppose that $Q$ is a scheme such that $v(Q)$ is maximal. Let $P$ be a node and assume that every block containing $P$ has size less than $t/2$. Then $v(P)<2\delta/t$.
\end{lemma}

\begin{proof2}
	Let the blocks of $Q$ that contain $P$ have sizes $r_1,...,r_u$. Then $\sum_i r_i=t+u-1$. Let $k$ be the minimal integer greater than 2 that is equal to some $r_i$ (or, if no such integer exists, then let $k$ be large enough that $\delta/k-\delta/(k+1)<\eta/2$). Let $R=\lfloor \frac{t-1}{2} \rfloor$. By assumption, $r_i\leq R$ for all $i$. Moreover, by the maximality of $v(Q)$ and Lemma \ref{app1}, we have the inequality $k\eta\geq \delta$ and therefore $\delta/k-\delta/(k+1)=\frac{\delta}{k(k+1)}\leq \frac{\eta}{k+1}<\eta/2$.
\medskip
	
\noindent \emph{Claim 1.} There exist positive integers $w$ and $q_1,...,q_w$ such that
	\begin{enumerate}[label=(\roman*)]
	\item $2\leq q_j\leq R$ for all $j$
	\item $\sum_j q_j=t+w-1$
	\item There is at most one $j$ for which $2<q_j<k$ and if there is any $i$ with $q_i=2$, then there is no $j$ with $2<q_j<k$.
	\item $v(P)\leq 1+\sum_j (\delta/q_j+(\alpha-1))$
	\item Either all but at most one $q_j$ are equal to $R$ or else $q_j\in \{2,R\}$ for all $j$
	\end{enumerate}
	
\noindent \emph{Proof of Claim 1.} Note that $v(P)=1+\sum_i (\delta/r_i+(\alpha-1))$. Define $w,q_1,q_2,...q_w$ to be the integers that maximize the quantity $1+\sum_j (\delta/q_j+(\alpha-1))$ subject to the conditions (i),(ii) and (iii). Since the $r_i$ satisfy (i),(ii),(iii), we get $v(P)\leq 1+\sum_j (\delta/q_j+(\alpha-1))$. We are left to prove (v), so let us suppose that it does not hold. There are two cases to consider.
\smallskip
	
	\emph{Case 1.} If there exists some $i$ with $q_i=2$, then there is a $j$ such that $q_j\not \in \{2,R\}$ and by (iii) we have $q_j\geq k$. Hence, $\delta/q_j-\delta/(q_j+1)<\eta/2$. After relabelling, we may assume that $j=w-1,i=w$. Now set $w'=w-1$, $q'_h=q_h$ for all $h\leq w-2$ and $q'_{w-1}=q_{w-1}+1$. Then $q'_1,...,q'_{w'}$ satisfy (i),(ii),(iii) and 
	\[1+\sum_{h\leq w} (\delta/q_h+(\alpha-1))<1+\sum_{h\leq w'} (\delta/q'_h+(\alpha-1)),\] 
which is a contradiction.
\smallskip	

	\emph{Case 2.} If there is no $i$ with $q_i=2$, then since (v) is assumed to fail, there must exist $i\neq j$ with $2<q_i\leq q_j<R$. Moreover, we may assume that $q_i$ is minimal among all $q_h$s. Without loss of generality, $i=w-1,j=w$. Now define $q'_h=q_h$ for all $h\leq w-2$, $q'_{w-1}=q_{w-1}-1$ and $q'_w=q_w+1$. Then $q'_1,...,q'_{w}$ satisfy (i),(ii),(iii) and 
	\[1+\sum_{h\leq w} (\delta/q_h+(\alpha-1))<1+\sum_{h\leq w} (\delta/q'_h+(\alpha-1)),\] 
which is a contradiction.
\medskip

	This completes the proof of Claim 1.
\medskip
	
\noindent \emph{Claim 2.} If $q_1,...,q_w$ satisfy the conditions (i),(ii),(v) in Claim 1, then 
\[1+\sum_{h\leq w} (\delta/q_h+(\alpha-1))< 2\delta/t.\]
	
\noindent \emph{Proof of Claim 2.} For $t\leq 13$, this is a straightforward check, which we performed using a computer program, since it would have taken inordinately long to do it by hand. (The code, written in Matlab, can be found at the end of the appendix.) So we shall assume that $t\geq 14$. Then $3R\geq 3\cdot\frac{t-2}{2}>t+2$, so there are at most two $q_j$s with $q_j=R$. Using (v), this leaves the following cases.
\smallskip
	
	Case 1: $q_j=2$ for all $j$
	
	Case 2: $q_1=R$	and $q_j=2$ for all $j\geq 2$
	
	Case 3a: $q_1=q_2=R=\frac{t-2}{2}$ and $q_3=q_4=q_5=2$ ($w=5$)
	
	Case 3b: $q_1=q_2=R=\frac{t-1}{2}$ and $q_3=q_4=2$	($w=4$)
	
	Case 4a: $q_1=q_2=R=\frac{t-2}{2},q_3=4$ ($w=3$)
	
	Case 4b: $q_1=q_2=R=\frac{t-1}{2},q_3=3$ ($w=3$)
	
\medskip
	
         By Lemmas \ref{app1} and \ref{app2} (d) we have $$(l-1)(\delta/2+(\alpha-1))<(\delta/l+(\alpha-1))$$ when $l=\frac{t-1}{2}$. Moreover, we have the inequality 
         \[\big(\delta/\Big(\frac{t-2}{2}\Big)+(\alpha-1)\big)+\frac{1}{2}\big(\delta/2+(\alpha-1)\big)<\big(\delta/\Big(\frac{t-1}{2}\Big)+(\alpha-1)\big),\] 
since this is equivalent to $\frac{2\delta}{(t-1)(t-2)}<\eta/4$, which holds because $(t-1)\eta\geq 2\delta$ and $t-2>4$. It is not hard to see that these two observations allow us to deduce all Cases 1-3 from Case 3b. To prove Case 3b, we need the inequality 
\[1+2(\delta/(\frac{t-1}{2})+(\alpha-1))-\eta< 2\delta/t,\] which is given in Lemma \ref{app2} (f).
	
	Clearly, Case 4a follows from Case 4b. To prove Case 4b, we need 
	\[1+2(\delta/(\frac{t-1}{2})+(\alpha-1))+(\delta/3+(\alpha-1))< 2\delta/t.\] 
Using $\alpha<1/2$ and $\delta<1$, it suffices to prove that $4/(t-1)-2/t\leq 1/6$, which holds for $t\geq 14$.	
\medskip

This completes the proof of Claim 2, and the two claims imply the lemma.
\end{proof2}

\begin{lemma} \label{contains}
	Suppose that $Q$ is a scheme such that its $v(Q)$ is as large as possible and such that the largest block $D$ of $Q$ has size at least $t/2$. Then $D$ has size $s$.
\end{lemma}

\begin{proof2}
	Suppose not. Pick a node $P$ with $P\not \in D$. Let $D$ have size $k\geq t/2$. Suppose that $P$ is contained in exactly $r$ large blocks. Define a scheme $Q'$ as follows. $Q'$ has the same blocks as $Q$ except that
	\begin{itemize}
		\item $P$ is removed from all large blocks,
		\item all small blocks containing $P$ and a node in $D$ are deleted,
		\item $P$ is added to $D$,
		\item the missing edges are now provided by distinct colours.
	\end{itemize}
	
	We now compare the values $v(Q)$ and $v(Q')$. The node $P$ is in only one large block in $Q'$ while it is in $r$ large blocks in $Q$. The number of small blocks containing $P$ is precisely $t-k-1$ in $Q'$ while it is at least $k-r$ in $Q$. That is because any large block containing $P$ contains at most one element of $D$. So 
	\begin{align*}v(Q')-v(Q)&\geq (r-1)(1-\alpha)+((t-k-1)-(k-r))(\delta+2(\alpha-1))\\
	&=(r-1)(1-\alpha)-(t-2k+(r-1))\eta\geq (r-1)(1-\alpha-\eta).\\
	\end{align*} 
But $1-\alpha-\eta=\delta-(1-\alpha)=1/2+(2s-1)\epsilon-(1/2+\epsilon)>0$. This contradicts the maximality of $v(Q)$ if $r\geq 2$.
	
	If $r=1$, then let the unique large block containing $P$ have size $l$. By assumption, $l\leq k$. Hence,  
	\[v(Q')-v(Q)=((t-k-1)-(t-l))(\delta+2(\alpha-1))=(k+1-l)\eta>0,\] 
a contradiction.
	
	If $r=0$, then 
	\[v(Q')-v(Q)=-(1-\alpha)-k(\delta+2(\alpha-1))=k\eta-(1-\alpha)\geq \frac{t}{2}\eta-(1-\alpha).\] 
But by Lemma \ref{app2} (d), this is at least $\delta-(1-\alpha)>0$. This is a contradiction and the lemma is proved.
\end{proof2}

\vspace{2mm}

We are ready to complete the proof of Lemma \ref{negscheme}.

\begin{proof}[\textnormal{\textbf{Proof of Lemma \ref{negscheme}.}}]
	We may assume that $v(Q)$ is maximal possible among all schemes $Q$. If $Q$ has a block of size $s$, then we are done by Lemma \ref{copycase}. Otherwise, by Lemma \ref{contains}, there is no block of size greater than or equal to $t/2$. But then Lemma \ref{locsum} and Lemma \ref{localneg} together imply that $v(Q)\leq t\frac{2\delta}{t}-(\delta+2\alpha)=\delta-2\alpha<0$.
\end{proof}

\appendix

\section{Appendix}

\begin{lemma} \label{app1}
	For any $k>2$, we have \begin{equation*}
	{k \choose 2}(\delta+2(\alpha-1))>\delta+k(\alpha-1) \Longleftrightarrow k\eta<\delta
	\end{equation*}
\end{lemma}

\begin{proof2}
\begin{align*}
{k \choose 2}(\delta+2(\alpha-1))&>\delta+k(\alpha-1)\\
\iff &(k-1)(\delta+2(\alpha-1))>2\delta/k+2(\alpha-1)\\
\iff &(k-1)\eta<2(1-\alpha)-2\delta/k=\eta+\delta(1-2/k)\\
\iff &(k-2)\eta<\delta(k-2)/k\\
\iff & k\eta<\delta.\\
\end{align*}
\end{proof2}

\begin{lemma} \label{app2}
	(a) $(t-s+1)\eta<\delta$ if and only if $(s,t)$ is regular.
	
	(b) $(t-s)\eta<\delta$ unless $t=2s-1$
	
	(c) $(t-s-1)\eta<\delta$
	
	(d) $(t-1)\eta> 2\delta$.
	
	(e) $\delta>2/3$.
	
	(f) $1+2(\delta/(\frac{t-1}{2})+(\alpha-1))-\eta< 2\delta/t$
\end{lemma}

\begin{proof2}
	Assume first that $(s,t)$ is regular. Then after some tedious calculations, one finds that (a) is equivalent to the inequality
	\begin{equation*}
	(s-2)(t-s-2)(2s-t-3)-t-3s+8>0
	\end{equation*}
The left hand side is a quadratic in $t$ with negative leading coefficient so it is enough to check that the inequality holds when $t=s+3$ and when $t=2s-4$.
	
	For $t=s+3$ we require $(s-2)(s-6)-4s+5>0$, which holds for $s\geq 11$, and for $t=2s-4$ we require $(s-2)(s-6)-5s+12>0$, which holds for $s\geq 11$. It therefore suffices to check the inequality for the pairs $(s,t)=(10,14)$ and $(s,t)=(10,15)$. This can be done by direct substitution. So (a) is proved (when $(s,t)$ is regular) which immediately implies (b) and (c).
	
	Now let us assume that $(s,t)$ is exceptional. Then the inequality $(t-s+c)\eta<\delta$ is equivalent to the inequality
	\begin{equation} \label{ceqn}
	(s-2)(t-s-c-1)(2s-t-2c-1)+(-2c^2-2c+1)s-t+4c^2+3c+1>0
	\end{equation}
When $c=1$, this says that $(s-2)(t-s-2)(2s-t-3)-3s-t+8>0$, so in order to prove (a) we need to show that this does not hold. For $t\in \{s+2,2s-3,2s-2,2s-1\}$ that is clear, since $(s-2)(t-s-2)(2s-t-3)\leq 0$. We are left to check that the inequality fails for the pairs $(7,10)$, $(8,11)$, $(8,12)$, $(9,12)$, $(9,13)$, $(9,14)$, $(10,13)$, and $(10,16)$. If $t=s+3$, then we need $s^2-12s+17\leq 0$ which indeed holds for $7\leq s\leq 10$. If $t=2s-4$, then we need $s^2-13s+24\leq 0$ which indeed holds for $7\leq s\leq 10$. We have only $(s,t)=(9,13)$ left to check. That is done by direct substitution.
	
	When $c=0$, then (\ref{ceqn}) says that $(s-2)(t-s-1)(2s-t-1)+s-t+1>0$. But if $s+2\leq t\leq 2s-2$, then the left hand side is minimal at $t=2s-2$ and there it takes value $(s-3)^2>0$. (Note that $s>3$ in this case.) This proves (b).
	
	When $c=-1$ in (\ref{ceqn}), then it says that $(s-2)(t-s)(2s-t+1)+s-t+2>0$. But the left hand side is minimal when $t=2s-1$, and then it is $2s^2-7s+7>0$. This proves (c).
	
	(d) In the regular case the statement is equivalent to the inequality
	\begin{equation*}
	2s^3-s^2t-5s^2+3st-t^2+s+3t-4> 0
	\end{equation*}
But in the regular case we have $2s-t\geq 4$, so $2s^3-s^2t\geq 4s^2$. Since $-s^2+3st-t^2\geq 0$ and $s+3t-4>0$, the statement follows.
	
	In the exceptional case, (d) is equivalent to the inequality
	\begin{equation*}
	(s-2)(2s-t)^2+(t-s-1)> 0
	\end{equation*}
which is clear.
	
	(e) In the regular case, the statement is equivalent to the inequality
	\[2s(t+s-5)(t-s-2)+2s^2-10s+6t-8>0,\] 
	which is easily seen to hold.
	
	In the exceptional case, it is equivalent to the inequality 
	\[(2s^2-5s)(t-s-2)+s^2-5s+2t-3>0,\] 
	which again clearly holds.
	
	(f) Since (by (d)) we have $\delta/(\frac{t-1}{2})<\eta$, this inequality reduces to
	
	$$ 1+2(\alpha-1)+2\delta/(t-1)<2\delta/t, $$
	
	or, equivalently, to
	
	$$ 2\alpha<1-\frac{2\delta}{t(t-1)}. $$ Expressing $\alpha$ in terms of $\delta$ and performing some routine algebraic manipulations, we find that we need to prove that
	
	$$ \frac{2(2s-1)}{t(t-1)}\delta<2\delta-1. $$
	
	Since $\delta<1$, the left hand side of this inequality is less than $4/t<1/3$ while the right hand side is greater than $1/3$, by part (e), so the proof is complete.
\end{proof2}

\medskip

Below we present the Matlab code that we used to perform the case check in the proof of Lemma \ref{localneg}.

\begin{verbatim}
% go through all pairs (s,t)
for t=5:13
    for s=(floor(t/2)+1):(t-2) 
        % these pairs are all exceptional
        alpha=((s-2)*(t-s)*(s-1)+s-1)/((2*s-3)*(t-s)*(s-1)+2*s-t); 
        delta=s-(2*s-1)*alpha;
        eta=2*(1-alpha)-delta;
        % bad will be changed to 1 if the inequality that we want 
        % to prove fails
        bad=0;     
        R=floor((t-1)/2);
        % j will count the number of q_h which are equal to R
        for j=0:4           
            a=(t-1)-j*(R-1);
            if 0<=a
                % in the following case every q_h is 2 or R
                v=1+a*(delta/2+alpha-1)+j*(delta/R+alpha-1);
                % check that our inequality holds with a suitably large
                % difference which can't be due to rounding errors
                if v>2*delta/(t)-10^(-3)                  
                    bad=1;
                end
            end
            if (2<=a+1) && (a+1<=R)
                % in the following case there is only one q_h that 
                % is not equal to R
                v=1+(delta/(a+1)+alpha-1)+j*(delta/R+alpha-1);  
                if v>2*delta/(t)-10^(-3)
                    bad=1;
                end
            end        
        end
        % tabulate the result: for each pair (s,t) we print 
        % whether the inequality failed (1) or not (0)
        fprintf('%5d %5d %5d \n',s,t,bad)
    end
end
\end{verbatim}


\section*{Acknowledgments} 
The second author was in Paris, supported by the Fondation Sciences Math\'ematiques de Paris, while this work was carried out. The first author held an FSMP Chair for the academic year 2017-8, also while this work was carried out, and would like to thank the FSMP for its support. He would also like to thank the \'Equipe d'Analyse Fonctionelle at Sorbonne Universit\'e for hosting him while he was in Paris. We are grateful to the two anonymous referees for their very careful reviews.



\begin{aicauthors}
\begin{authorinfo}[tim]
  Timothy Gowers\\
  Department of Pure Mathematics and Mathematical Statistics\\
  University of Cambridge, UK\\
  wtg10\imageat{}dpmms\imagedot{}cam\imagedot{}ac\imagedot{uk} \\
  \url{https://www.dpmms.cam.ac.uk/person/wtg10}
\end{authorinfo}
\begin{authorinfo}[oliver]
  Oliver Janzer\\
  Department of Pure Mathematics and Mathematical Statistics\\
  University of Cambridge, UK\\
  oj224\imageat{}cam\imagedot{}ac\imagedot{}uk\\
  \url{https://www.maths.cam.ac.uk/person/oj224}
\end{authorinfo}
\end{aicauthors}

\end{document}